\documentclass{article}
\usepackage{graphicx}
\usepackage{amsmath}
\usepackage{amsfonts}
\usepackage{amssymb}
\newtheorem{theorem}{Theorem}

\newtheorem{definition}[theorem]{Definition}

\begin{document}

\title{Geometric Aspects of Mirror Symmetry \\(with SYZ for Rigid CY manifolds)}
\author{Naichung Conan Leung\\School of Mathematics, \\University of Minnesota, \\Minneapolis, MN 55455, U.S.A..}
\maketitle
\begin{abstract}
In this article we discuss the geometry of moduli spaces of (1) flat bundles
over special Lagrangian submanifolds and (2) deformed Hermitian-Yang-Mills
bundles over complex submanifolds in Calabi-Yau manifolds.

These moduli spaces reflect the geometry of the Calabi-Yau itself like a
mirror. Strominger, Yau and Zaslow conjecture that the mirror Calabi-Yau
manifold is such a moduli space and they argue that the mirror symmetry
duality is a Fourier-Mukai transformation. We review various aspects of the
mirror symmetry conjecture and discuss a geometric approach in proving it.

The existence of rigid Calabi-Yau manifolds poses a serious challenge to the
conjecture. The proposed mirror partners for them are higher dimensional
generalized Calabi-Yau manifolds. For example, the mirror partner for a
certain K3 surface is a cubic fourfold and its Fano variety of lines is
birational to the Hilbert scheme of two points on the K3. This hyperk\"{a}hler
manifold can be interpreted as the SYZ mirror of the K3 by considering
\textit{singular} special Lagrangian tori.

We also compare the geometries between a CY and its associated generalized CY.
In particular we present a new construction of Lagrangian submanifolds.
\end{abstract}

\bigskip

\pagebreak 

\section{Introduction}

The mirror symmetry conjecture was proposed more than ten years ago from
string theorists. From a mathematical viewpoint, it is an amazing conjecture
about the geometry of Calabi-Yau manifolds.

Complex K\"{a}hler geometry studies manifolds $M$ with $U\left(  n\right)  $
holonomy. Using the complex structure $J$ and the Riemannian metric $g$, we
define a parallel non-degenerate two form $\omega,$
\[
\omega\left(  X,Y\right)  =g\left(  JX,Y\right)  .
\]
This is called the K\"{a}hler form and it defines a symplectic structure on
$M$.

In complex geometry we study objects such as complex submanifolds and
holomorphic vector bundles, or more generally coherent sheaves. Riemannian
metric on $M$ is then used to rigidify these objects. By the Wirtinger's
theorem, complex submanifolds are already rigidified, or calibrated. In
particular, they are absolute volume minimizers. To rigidify a holomorphic
vector bundle in geometry, we often look for a Hermitian-Yang-Mills
connection, namely its curvature tensor satisfies the equation
\[
F\wedge\omega^{n-1}=0\text{.}%
\]
A famous theorem of Donaldson, Uhlenbeck and Yau says that on a Mumford stable
bundle, there exists a unique Hermitian-Yang-Mills connection. In geometric
invariant theory, we need Gieseker stable bundle to construct algebraic moduli
space. On such bundles we have connections which satisfy a deformation of the
Hermitian-Yang-Mills equations \cite{L1}:
\[
\left[  e^{\frac{i}{2\pi}F+k\omega}Td\left(  M\right)  \right]  ^{2n}=C\left(
k\right)  \omega^{n}I_{E},
\]
for large enough $k$. We are going to absorb the factor $2\pi$ in the
definition of $F$. If the Todd class of $M$ is trivial and $k$ equals one,
then this equation becomes
\[
\left(  \omega+iF\right)  ^{n}=\text{constant.}%
\]
Note that both $\omega$ and $iF$ are real forms. From string theoretical
considerations, the following deformed Hermitian-Yang-Mills equation is
introduced in \cite{MMMS} in order to preserve supersymmetry:%

\[
\operatorname{Im}\left(  \omega+F\right)  ^{n}=0\text{.}%
\]
A supersymmetry B-cycle is defined to be a pair $\left(  C,E\right)  $ where
$C $ is a complex submanifold of $M$ and $E$ is a deformed
Hermitian-Yang-Mills bundle over $C$. In section \ref{GeomCycle} and the
appendix we will discuss the deformation theory of these B-cycles and
introduce the correlation function which is the following n-form on their
moduli space
\[
_{B}\Omega=\int_{C}Tr\mathbb{F}^{m}\wedge ev^{\ast}\Omega.
\]

On the other hand we can also study the symplectic geometry of $M$ using its
K\"{a}hler form $\omega$ as the symplectic form. Natural geometric objects are
Lagrangian submanifolds and their unitary flat bundles. To rigidify Lagrangian
submanifolds we need a parallel n-form. This requires $M$ to be a Calabi-Yau
manifold, namely its Riemannian metric has holonomy group in $SU\left(
n\right)  $. We denote its parallel holomorphic volume form as $\Omega$, it
satisfies
\[
\Omega\bar{\Omega}=\left(  -1\right)  ^{n\left(  n+1\right)  /2}2^{n}%
i^{n}\frac{\omega^{n}}{n!}.
\]
\newline This equation is equivalent to the Ricci flat condition for
Calabi-Yau manifolds. Now we define a supersymmetry A-cycle as a pair $\left(
C,E\right)  $ where $C$ is a special Lagrangian submanifold in $M$ and $E$ is
unitary flat bundle over $C$. A Lagrangian submanifold $C$ is called special
\cite{HL} if it is calibrated by $\operatorname{Re}\Omega$, namely
\[
\operatorname{Re}\Omega|_{C}=vol_{C}\text{.}%
\]
This is also equivalent to
\[
\operatorname{Im}\Omega|_{C}=0\text{.}%
\]
Special Lagrangian submanifolds minimize volume functional among submanifolds
representing the same homology class. In section \ref{GeomCycle}, we discuss
their deformation theory (after \cite{Mc}) and introduce the correlation
function which is also an n-form on the moduli space, namely%

\[
_{A}\Omega=\int_{C}Tr\left(  ev^{\ast}\omega+\mathbb{F}\right)  ^{n}\text{.}%
\]
On the symplectic side we have the Gromov-Witten theory which, roughly
speaking, counts the number of holomorphic curves of various genera in $M$. In
fact one should allow holomorphic curves to have real boundary whose images
lie on $C$ (see for example \cite{FOOO}). We call such holomorphic curves instantons.

Mirror symmetry conjectures predict that the complex geometry and the
symplectic geometry on Calabi-Yau manifolds are essentially equivalent. More
precisely for $M$ in a rather large class of Calabi-Yau manifolds, there
should exist another Calabi-Yau manifold $W$ in the same class called the
mirror manifold, this relation is reflexive; The complex geometry of $M$
should be equivalent to the symplectic geometry of $W$, with suitable
instanton corrections. Conversely the symplectic geometry of $M$ should be
equivalent to the complex geometry of $W$.

To understand \textit{why and how }these two kinds of geometry are
interchanged between mirror manifolds, we should first look at the semi-flat
case. The importance of the semi-flat case is first brought up by Strominger,
Yau and Zaslow in their paper \cite{SYZ}, which explains the mirror symmetry
from a physical/geometric viewpoint. This is now called the SYZ mirror
conjecture. The semi-flat case is then studied by Hitchin \cite{H1}, Gross
\cite{Gr2}, Yau, Zaslow and the author in \cite{LYZ} and also \cite{L3}. The
main advantage here is the absent of instantons:

We start with an affine manifold $D$, for simplicity we assume $D$ is just a
domain in $\mathbb{R}^{n}$. Let $\phi$ be a solution to the real
Monge-Amp\`{e}re equation
\begin{align*}
\det\left(  \nabla^{2}\phi\right)   &  =1,\\
\nabla^{2}\phi &  >0.
\end{align*}
Here $\nabla^{2}\phi$ is the Hessian of $\phi$ and it defines a Riemannian
metric on $D$ which we call a Cheng-Yau manifold because of their fundamental
result on the existence of such structure. Any such solution determines two
open Calabi-Yau manifolds, $TD$ and $T^{\ast}D$. Notice that $T^{\ast}D$
carries a canonical symplectic structure and $TD$ carries a canonical complex
structure because $D$ is affine. We can also compactify the fiber directions
by quotienting $TD$ and $TD^{\ast}$ with a lattice $\Lambda$ in $\mathbb{R}%
^{n}$ and its dual lattice $\Lambda^{\ast}$ in $\mathbb{R}^{n\ast}$
respectively and obtains mirror manifolds $M$ and $W$. The natural fibrations
on $M$ and $W$ over $D$ are both special Lagrangian fibrations.

The mirror transformation from $M$ to $W,$ and vice versa, is a generalization
of (i) the Fourier transformation on fibers of $M\rightarrow D $ together with
(ii) the Legendre transformation on the base $D$. In section \ref{MSsemiflat}
we will explain how the mirror transformation exchanges the complex geometry
and the symplectic geometry between $M$ and $W $.

The moduli space of complexified symplectic structures on $M$ is canonically
identified with the moduli space of complex structures on $W$; moreover this
map identifies various geometric structures on these moduli spaces, including
the two Yukawa couplings. Moduli spaces of certain supersymmetric A-cycles on
$M$ are canonically identified with moduli spaces of supersymmetric B-cycles
on $W$, moreover this map identifies various geometric structures on these
moduli spaces, including the two correlation functions. Holomorphic
automorphisms of $M$ are transformed to symplectic automorphisms of $W$ when
they are linear along fibers. On $M$ (and also $W$) there is an $\mathbf{sl}%
\left(  2\right)  \times\mathbf{sl}\left(  2\right)  $ action on its
cohomology groups, these two $\mathbf{sl}\left(  2\right)  $ actions are
interchanged under the mirror transformation from $M$ to $W$.

For general Calabi-Yau mirror manifolds $M$ and $W,$ it is conjectured that
above relationships should continue to hold after including instanton
corrections on the symplectic side. In section \ref{MirrorConj} we explain
various aspects of mirror symmetry conjectures. We also include a discussion
on a closely related conjecture of Gopakumar and Vafa.

Even though the Fourier/Legendre transformation approach provides a good
conceptural explanation to the mirror symmetry conjecture for Calabi-Yau
\textit{hypersurfaces} in toric varieties, there is a big puzzle when a
Calabi-Yau manifold is only a complete intersection in a toric variety because
such a Calabi-Yau may be a \textit{rigid} manifold, namely $H^{1}\left(
M,T_{M}\right)  =0$. Its mirror, if exists, would have $H^{1}\left(
W,T_{W}^{\ast}\right)  =0$ and therefore $W$ cannot be K\"{a}hler, a
contradiction! Nevertheless physicists still predict $M$ to have mirror
duality and its mirror partner is a higher dimensional Fano manifold of a
special kind.

In section \ref{RigidCY}, we provide a possible explanation of this mirror
symmetry conjecture for (possibly rigid) Calabi-Yau manifolds using the SYZ
approach. The basic idea is the complete intersection Calabi-Yau manifold
should be interpreted as a moduli space of B-cycles in the higher dimensional
Fano manifold. For example when $M$ is a Kummer K3 surface associated to the
product of two elliptic curves with complex multiplication, its physical
mirror partner would be the Fermat cubic fourfold $\bar{W}$ in $\mathbb{CP}%
^{5}$. It is a classical fact that the space of lines in $\bar{W}$, called the
Fano variety of lines, is birational to the Hilbert scheme of two points on
$M$. This hyperk\"{a}hler manifold $W$ could be interpreted as the SYZ mirror
to $M$. It is because each fiber in the natural special Lagrangian fibration
on the four torus becomes a singular torus with a node in $M$, a generic
deformation of it in $M$ is a smooth surface of genus two! The corresponding
moduli space of A-cycles in $M $ is the hyperk\"{a}hler manifold $W$.

It is a very important and interesting question to determine the complex and
symplectic geometry (or category) of the Calabi-Yau manifold $M$ from the one
on its associated higher dimensional Fano manifold $\bar{M}$. For example we
will present a new construction of a Lagrangian submanifold in $\bar{M}$ from
a (lower dimensional) Lagrangian submanifold in $M$. We would like to
understand how the Floer homology theory of Lagrangian intersections behaves
under this construction.

We should remark that the mirror symmetry for three dimensional Calabi-Yau
manifolds has much richer structure than their higher dimensional
counterparts. For example a conjecture of Gopakumar and Vafa \cite{GV2}
computes higher genus Gromov-Witten invariants of Calabi-Yau threefolds and
another conjecture of Ooguri and Vafa \cite{OV2} relates the large N
Chern-Simons theory of knots to counting holomorphic disks on a local
Calabi-Yau threefold. Certain parts of these conjectures are verfied
mathematically by Bryan and Pandharipande \cite{Pa}, \cite{BP}, Katz and Liu
\cite{KL}, Li and Song \cite{LS}. The basic reason behind these is the $G_{2}
$-structure on the seven dimensional manifold $M\times S^{1}$ (see e.g.
\cite{Lee Leung} and \cite{Leung intersect G2} and references therein).

\bigskip

\textit{Acknowledgments: The author thanks many people for valuable
discussions. Among them are J. Bryan, R. Donagi, B. Hassett, S. Katz, P.
Seidel, R. Thomas, C. Vafa, X.W. Wang, S.T. Yau, E. Zaslow and many more. The
paper is prepared when the author visited the Natural Center of Theoretical
Science, Tsing-Hua University, Taiwan in the summer of 2000. The author thanks
the center for providing an excellent research environment and support. This
project is also partially supported by a NSF grant, DMS-9803616.}

\textit{\pagebreak }

\section{\label{CYGeom}Geometry of Calabi-Yau manifolds}

Recall that a $2n$ dimensional compact Riemannian manifold $\left(
M,g\right)  $ with holonomy group inside $U\left(  n\right)  $ is a K\"{a}hler
manifold. It has a parallel complex structure $J$ and a parallel symplectic
structure $\omega$. By the classification result of Berger, possible special
holonomy groups inside $U\left(  n\right)  $ include $SU\left(  n\right)  $
and $Sp\left(  n/2\right)  $. Corresponding manifolds are called Calabi-Yau
manifolds and hyperk\"{a}hler manifolds. They can be characterized by the
existence of a parallel holomorphic volume form $\Omega$ and a parallel
holomorphic symplectic form $\eta$ respectively. The following table gives a
comparison among these geometries.%

\[%
\begin{tabular}
[c]{|l|l|l|l|}\hline
& K\"{a}hler & Calabi-Yau & hyperk\"{a}hler$%
\begin{array}
[c]{l}%
\,\\
\,
\end{array}
$\\\hline
Holonomy & $U\left(  n\right)  $ & $SU\left(  n\right)  $ & $Sp\left(
n/2\right)
\begin{array}
[c]{l}%
\,\\
\,
\end{array}
$\\\hline
Parallel forms & $\omega\in\Omega^{1,1}$ & $\omega\in\Omega^{1,1},\Omega
\in\Omega^{n,0}$ & $\omega\in\Omega^{1,1},\eta\in\Omega^{2,0}
\begin{array}
[c]{l}%
\,\\
\,
\end{array}
$\\\hline
Geometry$%
\begin{array}
[c]{l}%
\,\\
\,
\end{array}
$ & complex & cpx and sympl & $S^{2}$ family of cpx/sympl\\\hline
Action on $H^{\ast}\left(  M\right)  $ & \textbf{\ }$\mathbf{so}\left(
3\right)  $ & \textbf{\ }$\mathbf{so}\left(  4\right)  $ & $\mathbf{so}\left(
5\right)
\begin{array}
[c]{l}%
\,\\
\,
\end{array}
$\\\hline
\end{tabular}
\]
For the last row of the above table, the $\mathbf{so}\left(  3\right)  $
($=\mathbf{sl}\left(  2\right)  $) action on $H^{\ast}\left(  M\right)  $ is
the hard Lefschetz theorem. Notice that $\mathbf{so}\left(  4\right)
=\mathbf{sl}\left(  2\right)  \times\mathbf{sl}\left(  2\right)  $ and one of
these two $\mathbf{sl}\left(  2\right)  $ actions in the Calabi-Yau case come
from the K\"{a}hler geometry of $M$ as above. The existence of the other
$\mathbf{sl}\left(  2\right)  $ action is conjectural, this will be explained
in section \ref{MirrorConj}. In the semi-flat case this second $\mathbf{sl}%
\left(  2\right)  $ action arises from a variation of complex structures
toward the large complex limit point and it is described explicitly in
\cite{L3}.

\bigskip

\textbf{Yau's theorem}

In general it is difficult to find manifolds with special holonomy other than
$U\left(  n\right)  $. However in the $SU\left(  n\right)  $ case, we have the
celebrated theorem of Yau: Any compact K\"{a}hler manifold with trivial
canonical line bundle $K_{M}=\Lambda^{n}T_{M}^{\ast}$ admits a unique
K\"{a}hler metric with $SU\left(  n\right)  $ holonomy inside any given
K\"{a}hler class. Such a metric is called a Calabi-Yau metric.

\bigskip

\textbf{The complex Monge-Amp\`{e}re equation}

This theorem is obtained by solving a fully nonlinear elliptic equation on
$M\;$with $c_{1}\left(  M\right)  =0$: Suppose $\omega$ is any K\"{a}hler form
on $M$, then its Ricci form can be expressed as $Rc\left(  \omega\right)
=i\partial\bar{\partial}\log\omega^{n}$ with respect to any local holomorphic
coordinates on $\dot{M}.$ Since $c_{1}\left(  M\right)  =0$, we must have
$Rc\left(  \omega\right)  =i\partial\bar{\partial}F$ for some real valued
function $F$ on $M$. By the $\partial\bar{\partial}$-lemma, any other
K\"{a}hler form on $M$ in the same class must be of the form $\omega
+i\partial\bar{\partial}f$ for some real valued function $f$ on $M$. If this
new K\"{a}hler form has zero Ricci curvature, $Rc\left(  \omega+i\partial
\bar{\partial}f\right)  =0$, then
\begin{align*}
0  &  =i\partial\bar{\partial}\log\left(  \omega+i\partial\bar{\partial
}f\right)  ^{n}\\
&  =i\partial\bar{\partial}\log\frac{\left(  \omega+i\partial\bar{\partial
}f\right)  ^{n}}{\omega^{n}}+i\partial\bar{\partial}F\text{.}%
\end{align*}
On a closed manifold $M$, this implies that $\log\frac{\left(  \omega
+i\partial\bar{\partial}f\right)  ^{n}}{\omega^{n}}+F=C$ for some constant
$C$, which can be absorbed into $F$. Therefore the K\"{a}hler-Einstein
equation becomes the following complex Monge-Amp\`{e}re equation,
\[
\left(  \omega+i\partial\bar{\partial}f\right)  ^{n}=e^{-F}\omega^{n}\text{.}%
\]
By the work of Yau, this equation has a unique solution $f$, up to translation
by a constant. Equivalently, there is a unique Ricci flat metric on each
K\"{a}hler class if $c_{1}\left(  M\right)  =0$.

When the Ricci curvature is zero, $\int c_{2}\left(  M\right)  \omega^{n-2}$
becomes a positive multiple of the L$^{2}$ norm of the Riemannian curvature
tensor. Therefore if $M$ is a Calabi-Yau manifold, then
\[
\int c_{2}\left(  M\right)  \left[  \omega\right]  ^{n-2}\geq0,
\]
for any K\"{a}hler class $\left[  \omega\right]  $. Moreover it is zero if and
only if $M$ is covered by a flat torus.

\bigskip

$SU\left(  n\right)  $ \textbf{holonomy }

Notice that $c_{1}\left(  M\right)  =0$ implies the canonical line bundle
$K_{M}=\Lambda^{n}T_{M}^{\ast}$ of $M$ is topologically trivial. Suppose that
$K_{M}$ is holomorphically trivial, which is automatic if $M$ is simply
connected. Then there is a holomorphic section $\Omega$ of $K_{M}$. This is a
holomorphic $n$-form on $M$.

To see that Ricci flat metric on such manifold has holonomy inside $SU\left(
n\right)  $, we need the following result which is proven by using standard
Bochner type arguments: Every holomorphic $p$-form on a closed Ricci flat
K\"{a}hler manifold is parallel. In particular $\Omega$ is a parallel
holomorphic volume form. Using it, we can reduce the holonomy group from
$U\left(  n\right)  $ to $SU\left(  n\right)  $. We can also normalize
$\Omega$ so that
\[
\Omega\bar{\Omega}=\left(  -1\right)  ^{n\left(  n+1\right)  /2}2^{n}%
i^{n}\frac{\omega^{n}}{n!},
\]
where $\omega$ is the Ricci flat K\"{a}hler form on $M$. In fact any
K\"{a}hler form satisfying the above equation is automatically a Calabi-Yau
metric since $Rc=i\partial\bar{\partial}\log\omega^{n}$.

\bigskip

\textbf{Examples of Calabi-Yau manifolds}

Obvious examples of Calabi-Yau manifolds are complex tori $\mathbb{C}%
^{n}/\Lambda$ with the flat metrics, however their holonomy groups are
trivial. Smooth hypersurfaces in $\mathbb{CP}^{n+1}$ of degree $n+2$ have
trivial first Chern class, by Yau's theorem they are Calabi-Yau manifolds and
their holonomy groups equal $SU\left(  n\right)  $. When $n=2$, any smooth
quartic surface is actually hyperk\"{a}hler because $SU\left(  2\right)
=Sp\left(  1\right)  $. This is a K3 surface, namely a simply connected
K\"{a}hler surface with trivial first Chern class. It is known that all K3
surfaces form a connected (non-Hausdorff) moduli space. Moreover every
Calabi-Yau surface is either a K3 surface or a complex torus. For $n=3$ we
have the quintic threefold in $\mathbb{CP}^{4}$ and it is the most studied
Calabi-Yau threefold in mirror symmetry.

Among all Calabi-Yau hypersurfaces $M$ in $\mathbb{CP}^{n+1}$, there is a
natural choice. Namely $M$ is the union of all coordinate hyperplanes in
$\mathbb{CP}^{n+1}$:
\[
M=\left\{  \left(  z^{0},z^{1},...,z^{n+1}\right)  \in\mathbb{CP}^{n+1}%
:z^{0}z^{1}\cdots z^{n+1}=0\right\}  .
\]
It is singular! In a sense it is the most singular hypersurface and it is
called the large complex structure limit by physicists. It can also be
characterized mathematically by its mixed Hodge structure having maximal
unipotent monodromy (see for example \cite{CK}). Despite its non-smoothness,
it is a semi-flat Calabi-Yau space. We will explain this and its importances
in mirror symmetry in the following sections.

\label{Toric Eg}We can replace the projective space $\mathbb{CP}^{n+1}$ by a
Fano toric variety to generate many more explicit examples of Calabi-Yau
manifolds: Recall that a $\left(  n+1\right)  $ dimensional toric variety
$X_{\Delta}$ is a K\"{a}hler manifold with an effective $\left(
\mathbb{C}^{\times}\right)  ^{n+1}$-action. The analog of union of coordinate
hyperplanes in $\mathbb{CP}^{n+1}$ would be the union of its lower dimensional
$\left(  \mathbb{C}^{\times}\right)  ^{n+1}$-orbits. Again this is a singular
Calabi-Yau space and it can be deformed to a smooth, or mildly singular,
Calabi-Yau manifold if $X_{\Delta}$ is Fano or equivalently its associated
polytope $\Delta$ is reflexive.

One of many ways to define the polytope $\Delta$ for $X_{\Delta}$ is to use
the moment map in symplectic geometry. If we restrict the $\left(
\mathbb{C}^{\times}\right)  ^{n+1}$-action to its imaginary part, $\left(
S^{1}\right)  ^{n+1}$, then it preserves a symplectic form on $X_{\Delta}$ and
it has a moment map $\mu$,
\[
\mu:X_{\Delta}\rightarrow\mathbb{R}^{n+1}.
\]
The image $\mu\left(  X_{\Delta}\right)  $ is the associated polytope $\Delta
$. The pre-image of $\mu$ at a point on a $k$ dimensional open facet of
$\Delta$ is a $k$ dimensional torus $\left(  S^{1}\right)  ^{k}$. Therefore
the above singular Calabi-Yau hypersurface in $X_{\Delta}$ is just $\mu
^{-1}\left(  \partial\Delta\right)  $ and it admits a fibration by $n$
dimensional tori. This torus fibration will play a crucial role in the SYZ
mirror conjecture.

This construction has a natural generalization to Calabi-Yau complete
intersections in toric varieties, especially if every hypersurface involved is
semi-ample. Otherwise the Calabi-Yau manifold we construct as complete
intersection may be rigid, see section \ref{RigidCY} for its importance in
mirror symmetry.

\bigskip

\textbf{Moduli space of Calabi-Yau manifolds and Tian-Todorov results}

Infinitesimal deformations of the complex structures on $M$ are parametrized
by $H^{1}\left(  M,T_{M}\right)  .$ After fixing a holomorphic volume form on
the Calabi-Yau manifold $M$, we obtain a natural identification
\[
H^{1}\left(  M,T_{M}\right)  =H^{1}\left(  M,\Omega^{n-1}\right)
=H^{n-1,1}\left(  M\right)  .
\]

Given an infinitesimal deformation of complex structures on $M$, there may be
obstruction for it to come from a honest deformation. Various levels of
obstructions all lie inside $H^{2}\left(  M,T_{M}\right)  $. Tian \cite{Ti}
and Todorov \cite{Ti} show that all these obstructions vanish for Calabi-Yau
manifolds. As a corollary, the moduli space of complex structures for
Calabi-Yau manifolds is always smooth with tangent space equals $H^{n-1,1}%
\left(  M\right)  $.

Taking $L^{2}$ inner product of harmonic representatives in $H^{n-1,1}\left(
M\right)  $ with respect to the Calabi-Yau metric determines a K\"{a}hler
metric on this moduli space and it is called the Weil-Petersson metric. Tian
and Todorov show that its K\"{a}hler potential can be expressed in term of
variation of complex structures on $M$ in a simple manner. This is later
interpreted as a special geometry on the moduli space by Strominger. It also
plays an important role in constructing a Frobenius structure on an extended
moduli space (see for example \cite{CK} for more details).

\bigskip

\textbf{The definition of }$B$\textbf{-fields}

A purpose of introducing $B$-fields is to complexify the space of symplectic
structures on $M$, the conjectural mirror object to the space of complex
structures on $W$ which has a natural a complex structure.

The usual definition of a $B$-field is a harmonic form of type $\left(
1,1\right)  $ with respect to the Calabi-Yau metric. The harmonicity condition
does not respect mirror symmetry unless we are in the large complex and
K\"{a}hler structure limit. In this limit, both the complex and real
polarizations are essentially equivalent. Away from this limit, we define the
B-field $\beta$ in two different ways depending one whether we use the real or
complex polarization. In any case, $\beta$ would be a closed two form of type
$\left(  1,1\right)  $. If we use the complex polarization, then we require
$\beta$ to satisfy the equation $\operatorname{Im}\left(  \omega
+i\beta\right)  ^{n}=0$. We can show that if $M$ satisfies the following
condition.
\begin{align*}
\omega^{n}  &  =i^{n}\Omega\bar{\Omega}\\
\operatorname{Im}e^{i\theta}\left(  \omega+i\beta\right)  ^{n}  &  =0\\
\operatorname{Im}e^{i\phi}\Omega &  =0\text{ on the zero section.}%
\end{align*}
Then in the semi-flat case, under the mirror transformation to $W$, they
becomes
\begin{align*}
\omega_{W}^{n}  &  =i^{n}\Omega_{W}\bar{\Omega}_{W}\\
\operatorname{Im}e^{i\theta}\Omega_{W}  &  =0\text{ on the zero section.}\\
\operatorname{Im}e^{i\phi}\left(  \omega_{W}+i\beta_{W}\right)  ^{n}  &  =0.
\end{align*}

When we expand the second equation for small $\beta$ and assume the phase
angle is zero, we have
\[
\left(  \omega+i\varepsilon\beta\right)  ^{n}=\omega^{n}+i\varepsilon
n\beta\omega^{n-1}+O\left(  \varepsilon^{2}\right)  \text{.}%
\]
If we linearize this equation, by deleting terms of order $\varepsilon^{2}$ or
higher, then it becomes the harmonic equation for $\beta$:
\[
\beta\omega^{n-1}=c^{\prime}\omega^{n}\text{.}%
\]
That is (1) $\omega$ is a Calabi-Yau K\"{a}hler form and (2) $\beta$ is a
harmonic real two form. This approximation is in fact the usual definition of
a B-field.

\bigskip

\textbf{A remark}

Here we explain why we usually consider the complexified K\"{a}hler class to
be a cohomology class in $H^{1,1}\left(  M,\mathbb{C}\right)  /H^{1,1}\left(
M,\mathbb{Z}\right)  $ rather than inside $H^{1,1}\left(  M,\mathbb{C}\right)
$. Every element in $H^{1,1}\left(  M,\mathbb{Z}\right)  =H^{1,1}\left(
M\right)  \cap H^{2}\left(  M,\mathbb{Z}\right)  $ is the first Chern class of
a holomorphic line bundle $L$ on $M$. As we will see in section
\ref{GeomCycle}, $L$ determines a supersymmetric $B$-cycle on $M$ if its
curvature tensor $F\in\Omega^{2}\left(  M,i\mathbb{R}\right)  $ satisfies the
equation
\[
\operatorname{Im}e^{i\theta}\left(  \omega+i\left[  \beta+\frac{i}{2\pi
}F\right]  \right)  ^{n}=0\text{,}%
\]
for some constant $C$. Therefore shifting the cohomology class of $-i\left(
\omega+i\beta\right)  $ by an element in $H^{1,1}\left(  M,\mathbb{Z}\right)
$ is roughly the same as tensoring a supersymmetric cycle by a holomorphic
line bundle. The correct way to look at this issue is to consider all B-cycles
whose cohomology classes $\left[  Tr\left(  \omega/i+\left[  \beta+\frac
{i}{2\pi}F\right]  \right)  \right]  $ are the same in $H^{1,1}\left(
M,\mathbb{C}\right)  $ modulo $H^{1,1}\left(  M,\mathbb{Z}\right)  $. In this
respect the complexified Calabi-Yau metric on $M$ is the same as the trivial
line bundle being a supersymmetric $B$-cycle in $M$. In fact more precisely a
complexified K\"{a}hler class would be in $H^{1,1}\left(  M,\mathbb{C}%
/\mathbb{Z}\right)  $.

\bigskip

Another approach is to the B-field is via the real polarization of $M$, namely
a special Lagrangian fibration on $M$. The distinction here is the complex
conjugation will be replaced by a real involution, by sending the fiber
directions to its negative, namely $dx^{j}\rightarrow dx^{j}$ and
$dy^{j}\rightarrow-dy^{j}$. In the large complex and K\"{a}hler structures
limit, these two involutions are the same. We denote the real involution of
$\Omega$ by $\widehat{\Omega}$. Suppose that $\omega$ is a K\"{a}hler form on
$M$ and $\beta$ is a closed real two from on $M$ of type $\left(  1,1\right)
$. Then $\omega^{\mathbb{C}}=\omega+i\beta$ is called the complexified
Calabi-Yau K\"{a}hler class of $M$ if $\left(  \omega^{\mathbb{C}}\right)
^{n} $ is a nonzero constant multiple of $i^{n}\Omega\wedge\widehat{\Omega}$.
\[
\left(  -1\right)  ^{n\left(  n+1\right)  /2}2^{n}i^{n}\left(  \omega
^{\mathbb{C}}\right)  ^{n}=\Omega\wedge\widehat{\Omega}\text{.}%
\]
We called this the complexified complex Monge-Amp\`{e}re equation. We will see
explicitly in the semi-flat case that we need these modified definitions of
B-fields to exchange complex structures and complexified K\"{a}hler structures
of $M$ and $W$ (see \cite{L3} for details).

\bigskip

\pagebreak 

\section{\label{GeomCycle}Supersymmetric cycles and their moduli}

In this section we study geometric objects in $M$ and their moduli spaces in
complex geometry and symplectic geometry.

In complex geometry we study geometric objects such as complex submanifolds
and holomorphic vector bundles. In symplectic geometry we study Lagrangian
submanifolds. It is also natural to include flat bundles on these Lagrangian
submanifolds as an analog of holomorphic bundles on complex submanifolds on
the complex side.

From string theory considerations, Marino, Minasian, Moore and Strominger
argue in \cite{MMMS} that supersymmetry imposes metric constraints on these
geometric objects. On the complex side, they need to satisfy a deformation of
Hermitian-Yang-Mills equations. Such deformations are similar to those studied
by the author in \cite{L1} and \cite{L2} related to Gieseker stability of
holomorphic vector bundles. On the symplectic side, Lagrangian submanifolds
that preserve supersymmetry would be calibrated, the so-called special
Lagrangian submanifolds. This type of submanifolds are introduced and studied
by Harvey and Lawson in \cite{HL} related to minimal submanifolds inside
manifolds with special holonomy. These two types of objects are called
B-cycles and A-cycles because of their relationships with B-model and A-model
in string theory.

\bigskip

\textbf{B-cycles: Hermitian-YM bundles over complex submanifolds}

On B-side we study complex geometry of $M$. We do not need the Calabi-Yau
assumption on $M$ until we define the correlation function on their moduli spaces.

\begin{definition}
Let $M$ be a K\"{a}hler manifold with complexified K\"{a}hler form
$\omega^{\mathbb{C}}$, we call a pair $\left(  C,E\right)  $ a supersymmetry
B-cycle, or simply B-cycle, if $C$ is a complex submanifold of $M$ of
dimension $m$. $E$ is a holomorphic vector bundle on $C$ with a Hermitian
metric whose curvature tensor $F$ satisfies the following \textit{deformed
Hermitian-Yang-Mills equations }on $C$:
\[
\operatorname{Im}e^{i\theta}\left(  \omega^{\mathbb{C}}+F\right)  ^{m}=0,
\]
for some constant angle $\theta$ which is called the \textit{phase angle}.
\end{definition}

This equation is introduced in \cite{MMMS}. Recall that $E$ being holomorphic
is equivalent to its curvature tensor satisfies the integrability condition
$F^{2,0}=0$. From the complex geometry point of view, it is natural to include
those pairs $\left(  C,E\right)  $'s which are singular, namely a coherent
sheaf on $M$. Nevertheless it is unclear how to impose the deformed
Hermitian-Yang-Mills equations on such singular objects.

If we replace $\omega^{\mathbb{C}}$ by a large multiple $N\omega^{\mathbb{C}}
$, then the leading order term for the deformed Hermitian-Yang-Mills equation
becomes
\[
\operatorname{Im}e^{i\theta}F\wedge\left(  \omega^{\mathbb{C}}\right)
^{m-1}=0\text{.}%
\]
If the B-field and the phase angle $\theta$ are both zero then this equation
becomes the Hermitian-Yang-Mills equation for the vector bundle $E$:%

\[
\Lambda F=0.
\]

By the theorem of Donaldson, Uhlenbeck and Yau, an irreducible holomorphic
bundle $E$ over a compact K\"{a}hler manifold $C$ admits a
Hermitian-Yang-Mills connection if and only if the bundle is Mumford stable.
Moreover such a connection is unique. We recall that a holomorphic bundle $E$
is called Mumford stable if for any proper coherent subsheaf $S$ of $E$, we
have
\[
\frac{1}{rank\left(  S\right)  }\int c_{1}\left(  S\right)  \omega^{n-1}%
<\frac{1}{rank\left(  E\right)  }\int c_{1}\left(  E\right)  \omega^{n-1}.
\]
By the arguments in \cite{L1}, $E$ admits a solution to $\operatorname{Im}%
\left(  k\omega+F\right)  ^{m}=0$ for large enough $k$ is given by the
following stability notion: For any proper coherent subsheaf $S$ of $E$, we
require
\[
\frac{1}{rk\left(  S\right)  }\operatorname{Im}\int Tr\left(  k\omega
+F_{S}\right)  ^{m}<\frac{1}{rk\left(  E\right)  }\operatorname{Im}\int
Tr\left(  k\omega+F_{E}\right)  ^{m}%
\]
for all large enough $k$.

The deformation theory of B-cycles and the geometry of their moduli space
$_{B}\mathcal{M}\left(  M\right)  $ will be discussed in the appendix. For
example the tangent space to $_{B}\mathcal{M}\left(  M\right)  $ is
parametrized by deformed harmonic forms and there is a correlation function
$_{B}\Omega$ on $_{B}\mathcal{M}\left(  M\right)  $ defined as
\[
_{B}\Omega=\int_{C}Tr\mathbb{F}^{m}\wedge ev^{\ast}\Omega.
\]
They will play an important role in the mirror symmetry conjecture.

\bigskip\qquad

\textbf{A-cycles: Flat bundles on SLag submanifolds }\qquad

Special Lagrangian submanifolds are introduced by Harvey and Lawson in
\cite{HL} as calibrated submanifolds. Their original paper is an excellent
reference for the subject. It turns out that such objects are supersymmetric
cycles when coupled with deformed flat connections (see \cite{BBS}).

\begin{definition}
Let $M$ be a Calabi-Yau manifold of dimension $n$ with complexified K\"{a}hler
form $\omega^{\mathbb{C}}$ and holomorphic volume form $\Omega$. We called a
pair $\left(  C,E\right)  $ a supersymmetry A-cycle, or simply
\textit{A-cycle}, if (i) $C$ is a special Lagrangian submanifold of $M$,
namely $C$ is a real submanifold of dimension $n$ with
\[
\omega|_{C}=0,
\]
and
\[
\operatorname{Im}e^{i\theta}\Omega|_{C}=0.
\]

(ii) $E$ is a unitary vector bundle on $C$ whose curvature tensor $F$
satisfies the deformed flat condition,
\[
\beta|_{C}+F=0\text{.}%
\]
\end{definition}

Note that the Lagrangian condition and deformed flat equation can be combined
into one complex equation on $C$:\footnote{Since $\beta$ may represent a
non-integral class, we need to use generalized connections on the bundle as
studied in \cite{L2}.}
\[
\omega^{\mathbb{C}}+F=0.
\]
These conditions are very similar to those defining a B-cycle. There is also a
symplectic reduction picture on the A-side as introduced and studied by
Donaldson in \cite{D} and Hitchin in \cite{H2}. Thomas has a different
symplectic reduction picture and he compared it with the B-side story in
\cite{Th1}. The role of various objects are reversed - a mirror phenomenon.

\bigskip

\textbf{Moduli space of A-cycles and their correlation functions}

If we vary an A-cycle $\left(  C,E\right)  $ while keeping $C$ fixed in $M$,
then the moduli space of these objects is isomorphic to the moduli space of
flat connections on $C$. For simplicity we assume the B-field vanishes.
However there is no canonical object analog to the trivial flat connection
when the B-field is nonzero. Then this moduli space is independent of how $C $
sits inside $M$ because of its topological nature.

So we reduce the problem to understanding deformations of special Lagrangian
submanifolds and this problem has been solved by McLean in \cite{Mc}. Using
the Lagrangian condition, a normal vector field of $C$ can be identified with
an one form on $C$. McLean shows that infinitesimal deformations of $C$ are
parametrized by the space of harmonic one forms on $C$, moreover, all higher
order obstructions to deformations vanish. In particular every infinitesimal
deformation comes from a honest family of special Lagrangian submanifolds in
$M$. By Hodge theory the tangent space of moduli space of A-cycles,
$_{A}\mathcal{M}\left(  M\right)  $, at $\left(  C,E\right)  $ is given by
$H^{1}\left(  C,\mathbb{R}\right)  \times H^{1}\left(  C,ad\left(  E\right)
\right)  $.

When $E$ is a line bundle, we have a natural isomorphism $H^{1}\left(
C,ad\left(  E\right)  \right)  =H^{1}\left(  C,i\mathbb{R}\right)  $ and this
tangent space is
\[
H^{1}\left(  C,\mathbb{R}\right)  +H^{1}\left(  C,i\mathbb{R}\right)
=H^{1}\left(  C,\mathbb{C}\right)  .
\]
Therefore $_{A}\mathcal{M}\left(  M\right)  $ has a natural almost complex
structure which is in fact integrable (see for example \cite{H1}). By the
symplectic reduction procedure mentioned earlier, it also carries a natural
symplectic structure.

The correlation function is a degree $n$ form on this moduli space. First we
define a degree-$n$ closed form on the whole configuration space $Map\left(
C,M\right)  \times\mathcal{A}_{C}\left(  E\right)  $ as follows,
\[
_{A}\Omega=\int_{C}Tr\left(  ev^{\ast}\omega+\mathbb{F}\right)  ^{n}\text{.}%
\]
Here $ev:C\times Map\left(  C,M\right)  \rightarrow M$ is the evaluation map
and $\mathbb{F}$ is the curvature of the universal connection on
$C\times\mathcal{A}_{C}\left(  E\right)  $.

Since the tangent spaces of the moduli of special Lagrangians and the moduli
of flat connections can both be identified with the space of harmonic one
forms, a tangent vector of the moduli space $_{A}\mathcal{M}\left(  M\right)
$ is a complex harmonic one form $\eta+i\mu$. Then $_{A}\Omega$ can be written
explicitly as follows
\[
_{A}\Omega\left(  \eta_{1}+i\mu_{1},...,\eta_{m}+i\mu_{m}\right)  =\int
_{C}Tr\left(  \eta_{1}+i\mu_{1}\right)  \wedge...\wedge\left(  \eta_{m}%
+i\mu_{m}\right)  .
\]

\bigskip

\textbf{Instanton corrections}

As we have seen the A-side story is much easier to describe so far because the
equation involving $F$ is linear. However this correlation function is only
valid in the so-called classical limit. To realize the mirror symmetry, we
need to modify these structures by adding suitable contributions coming from
holomorphic disks whose boundaries lie on $C$. These are called instanton
corrections. In this paper we are not going to describe this aspect.

\bigskip

\textbf{A remark}

So far we have assumed $C$ and $E$ are always smooth objects. But if we want
to study the whole moduli space, it is natural to look for a suitable
compactification of it. In that case singular objects are unavoidable. On the
A-side, this problem is fundamental for understanding mirror symmetry via
special Lagrangian fibrations as proposed by Strominger, Yau and Zaslow
\cite{SYZ}. On the B-side, even when $C$ is a curve and $E$ is a line bundle,
understanding how to compactify this moduli space is crucial in the study of
the Gopakumar-Vafa conjecture \cite{GV2}.

\bigskip\newpage

\bigskip

\section{\label{MirrorConj}Mirror symmetry conjectures}

Roughly speaking the mirror symmetry conjecture says that for a Calabi-Yau
manifold $M$, there is another Calabi-Yau manifold $W$ of the same dimension
such that
\[
\fbox{$%
\begin{array}
[c]{ccc}%
\begin{array}
[c]{c}%
\text{Symplectic geometry}\\
\text{on }M
\end{array}
& \cong &
\begin{array}
[c]{c}%
\text{Complex Geometry}\\
\text{on }W
\end{array}
\end{array}
$}%
\]
and
\[
\fbox{$%
\begin{array}
[c]{ccc}%
\begin{array}
[c]{c}%
\text{Complex geometry}\\
\text{on }M
\end{array}
& \cong &
\begin{array}
[c]{c}%
\text{Symplectic Geometry}\\
\text{on }W.
\end{array}
\end{array}
$}%
\]
The conjecture cannot be true as stated because rigid Calabi-Yau manifolds are
counterexamples to it. From a string theory point of view, the conjecture
should hold true whenever the Calabi-Yau manifold admits a deformation to a
large complex structure limit. This includes Calabi-Yau hypersurfaces inside
Fano toric varieties. From a mathematical point of view, this should
correspond to the existence of a certain special Lagrangian fibration on $M$
as proposed in \cite{SYZ}. The conjecture might hold true for an even larger
class of Calabi-Yau manifolds. For example certain aspects of it is
conjectured to hold true even for rigid Calabi-Yau manifolds (see section
\ref{RigidCY}).

Other basic aspects of mirror symmetry are: (i) On the symplectic side, there
are instanton corrections. These are holomorphic curves or holomorphic disks
with boundaries on the special Lagrangian. (ii) On the complex side, the
structures are nonlinear. For example the deformed Hermitian-Yang-Mills
equations and the natural symplectic form on the moduli space of B-cycles are
nonlinear in the curvature tensor. Since one of the main purposes of this
paper is to understand how mirror symmetry works when there are no instantons
present, we will not spend much effort describing instanton effects even
though they are equally important in the whole subject.

In the early 90's, mirror symmetry mainly concern identification of the moduli
spaces of symplectic structures on $M$ and complex structures on $W$. The
information on the instanton corrections in this case predicts the genus zero
Gromov-Witten invariants of $M$.

Kontsevich \cite{K1} proposes a homological mirror symmetry conjecture. In
1996 Strominger, Yau and Zaslow \cite{SYZ} proposed a geometric mirror
symmetry conjecture (see also Vafa's paper \cite{V2}). These conjectures
relate the space of A-cycles in $M$ and the space of B-cycles in $W$. As an
analogy of topological theories, the homological conjecture is about singular
cohomology, or rational homotopy theory and the geometric conjecture is about
harmonic forms. Singular cohomology class can always be represented by a
unique harmonic form by Hodge theory. The corresponding results we need for
mirror symmetry would roughly be Donaldson, Uhlenbeck and Yau theorem and
Thomas conjecture (\cite{Th2}, \cite{TY}).

The key breakthrough is realizing $W$ as a moduli space of supersymmetric
A-cycles in $M$, as discovered in \cite{SYZ}. Before SYZ, many mirror
manifolds were constructed by Batyrev but this combinatorial construction does
not give us any insight into why the geometries of $M$ and $W$ are related at all.

From this geometric point of view, we can now understand why mirror symmetry
happens, at least when there are no instanton corrections (see section
\ref{MSsemiflat}). We can also enlarge the mirror symmetry conjecture to
include more geometric properties. For instance we will state a conjecture
regarding an $sl\left(  2\right)  \times sl\left(  2\right)  $ action on
cohomology groups of moduli spaces of supersymmetric cycles, which include in
particular $M$ and $W$.

Now we start from the beginning of mirror symmetry.

\bigskip

\textbf{The paper by Candelas, de la Ossa, Green and Parkes}

The mirror symmetry conjecture in mathematics begins with the fundamental
paper \cite{COGP} by Candelas, de la Ossa, Green and Parkes. They analyze the
mirror pair $M$ and $W$ constructed by Greene and Plesser. Here $M$ is a
Fermat quintic threefold in $\mathbb{P}^{4}$ and $W$ is a Calabi-Yau
resolution of a finite group quotient of $W$. In this paper they compare the
moduli space of complexified K\"{a}hler structures on $M$ and the moduli space
of complex structures on $W$. Both moduli spaces are one dimensional in this
case. The key structure on these moduli spaces is the Yukawa coupling. On the
complex side, the Yukawa coupling is determined by a variation of Hodge
structures on $W$ and it can be written down explicitly by solving a
hypergeometric differential equation. On the symplectic side, the Yukawa
coupling consists of two parts: classical and quantum. The classical part is
given simply by the cup product of the cohomology of $M$. The quantum part is
a generating function of genus zero Gromov-Witten invariants of $M$.

Before we can compare these two Yukawa couplings, we first need an explicit
map identifying the two moduli spaces. This mirror map cannot be arbitrary
because it has to transform naturally under the monodromy action. In this
case, this mirror map is determined in \cite{COGP}. Therefore the mirror
conjecture which identifies the two Yukawa couplings gives highly nontrivial
predictions on the genus zero Gromov-Witten invariant of $M$ by solving a
hypergeometric differential equation. The Gromov-Witten invariants of $M$ are
very difficult to determine even in low degree. At that time, only those of
degree not more than three were known mathematically. They are 2875, 609250,
317206375 etc.

\bigskip

\textbf{Gromov-Witten invariants}

Before stating the conjecture, let us briefly review Gromov-Witten theory.
Here we will only consider three dimensional Calabi-Yau manifolds $M$. Roughly
speaking the Gromov-Witten invariant $N_{C}^{g}\left(  M\right)  $, defined in
\cite{RT}, counts the number of genus $g$ curves in $M$ representing the
homology class $C\in H_{2}\left(  M,\mathbb{Z}\right)  $. To handle the
problem of possible non-reducedness and bad singularities of these curves, we
need to study instead holomorphic maps from genus $g$ stable curves to $M$.
The corresponding moduli space has expected dimension zero. If it is indeed a
finite number of smooth points, then its cardinality is $N_{C}^{g}\left(
M\right)  $. Otherwise one can still define a virtual fundamental class (see
\cite{LT} for example) and obtain a symplectic invariant $N_{C}^{g}\left(
M\right)  \in\mathbb{Q}$\footnote{These numbers can be fractional because of
automorphisms of stable maps.}. This is nonetheless a very nontrivial result.
On the other hand, computing these numbers $N_{C}^{g}\left(  M\right)  $ is a
very difficult mathematical problem.

Recently Gopakumar and Vafa \cite{GV2} proposed a completely different method
to obtain these $N_{C}^{g}\left(  M\right)  $ using the cohomology of the
moduli space of curves in the class $C\in H_{2}\left(  M,\mathbb{Z}\right)  $
together with flat $U\left(  1\right)  $ bundles over these curves.

Counting the number of curves in $M$ is an old subject in algebraic geometry,
called enumerative geometry. For example Clemens conjectured there are only
finite number of rational curves of any fixed degree on a generic quintic
threefold. However these $N_{C}^{0}\left(  M\right)  $'s do not really count
the number of rational curves in $M$. The same also applies to the higher
genus count too. For example, if there is only one curve in the class $C$ and
it is a smooth rational curve with normal bundle $O\left(  -1\right)  \oplus
O\left(  -1\right)  $. We have $N_{C}^{0}\left(  M\right)  =1$ as expected;
however, $N_{dC}^{0}\left(  M\right)  =\frac{1}{d^{3}}$ for every positive
integer $d$. This is the so-called multiple cover formula, conjectured in
\cite{COGP}, given a mathematical reasoning by Morrison and Plesser in
\cite{AM} and then proven rigorously by Voisin in \cite{Vo}.

So if all rational curves in $M$ are smooth and with normal bundle $O\left(
-1\right)  \oplus O\left(  -1\right)  $, then the genus zero Gromov-Witten
invariants do determine the number of rational curves in $M$. However even for
a generic quintic threefold in $\mathbb{P}^{4}$, this is not true. Vainsencher
\cite{Vai} shows that generic quintic threefold has degree five rational
curves with six nodes. There are in fact 17601000 such nodal curves. Therefore
to understand the enumerative meaning of the Gromov-Witten invariants, we need
know how each rational curve contributes. In \cite{BKL}, Bryan, Katz and the
author determine such contributions when the rational curve $C$ (i) has one
node and it is superrigid (analog of $O\left(  -1\right)  \oplus O\left(
-1\right)  $ for nodal curves) or (ii) is a smooth contractible curve. In the
second case, the normal bundle of $C$ must be either $O\oplus O\left(
-2\right)  $ or $O\left(  1\right)  \oplus O\left(  -3\right)  $.

\bigskip

\textbf{Candelas-de la Ossa-Green-Parkes mirror conjecture}

This famous conjecture has been discussed in details in many places (see e.g.
\cite{CK} for details) and our discussions will be very brief. Suppose that
$M$ and $W$ are mirror manifolds. Then the variation of complexified
symplectic structures of $M$ with instanton corrections should be equivalent
to the variation of complex structures of $W$. The infinitesimal deformation
spaces are given by $H^{1}\left(  M,T_{M}^{\ast}\right)  =H^{1,1}\left(
M\right)  $ and $H^{1}\left(  W,T_{W}\right)  =H^{n-1,1}\left(  W\right)  $
respectively. In particular we would have $\dim H^{1,1}\left(  M\right)  =\dim
H^{n-1,1}\left(  W\right)  $. More generally, it is conjectured that
\[
H^{q}\left(  M,\Lambda^{p}T_{M}^{\ast}\right)  =H^{q}\left(  W,\Lambda
^{p}T_{W}\right)  ,
\]
as they can be interpreted as tangent spaces of extended deformation problems
\cite{K1}. In particular we would have equalities for Hodge numbers%

\[
\dim H^{p,q}\left(  M\right)  =\dim H^{n-p,q}\left(  W\right)  .
\]
As a corollary, the conjecture implies equality of Euler characteristics up to
sign,
\[
\chi\left(  M\right)  =\left(  -1\right)  ^{n}\chi\left(  W\right)  \text{.}%
\]

The next step would be an identification of the moduli space of complexified
symplectic structures on $M$ and the moduli space of complex structures on $W
$, at least near a large complex structure limit point. This is called the
mirror map in the literature.

On each moduli space, there is degree $n$ form called the Yukawa coupling. The
mirror symmetry conjecture says that the Yukawa coupling on the moduli space
of complexified symplectic structures on $M$ is pulled back to the Yukawa
coupling on the moduli space of complex structures on $W$ by the mirror map.

On the $M$ side the classical part of Yukawa coupling is defined by the
natural product structure on cohomology groups:
\[%
\begin{array}
[c]{cc}%
_{A}\mathcal{Y}_{M}^{cl}: & \Lambda^{n}H^{1}\left(  M,T_{M}^{\ast}\right)
\rightarrow H^{n}\left(  M,\Lambda^{n}T_{M}^{\ast}\right)  =\mathbb{C}\\
& _{A}\mathcal{Y}_{M}^{cl}\left(  \zeta_{1},\zeta_{2},...,\zeta_{n}\right)
=\int_{M}\zeta_{1}\wedge\zeta_{2}\wedge...\wedge\zeta_{n}.
\end{array}
\]
With instanton corrections by genus zero Gromov-Witten invariants $N_{d}^{0}$,
the full Yukawa coupling when $n=3$ is%

\[
_{A}\mathcal{Y}_{M}^{cl}\left(  \zeta_{1},\zeta_{2},\zeta_{3}\right)
=\int_{M}\zeta_{1}\cup\zeta_{2}\cup\zeta_{3}+\sum_{d\in H^{2}\left(
M,\mathbb{Z}\right)  \backslash0}N_{d}^{0}\left(  \zeta_{1},\zeta_{2}%
,\zeta_{3}\right)  \frac{q^{d}}{1-q^{d}}.
\]

On the $W$ side, the Yukawa coupling is again defined by the natural product
structure on cohomology groups:
\[%
\begin{array}
[c]{cc}%
_{B}\mathcal{Y}_{W}: & \Lambda^{n}H^{1}\left(  W,T_{W}\right)  \rightarrow
H^{n}\left(  W,\Lambda^{n}T_{W}\right)  =\mathbb{C}\\
& _{B}\mathcal{Y}_{W}\left(  \zeta_{1},\zeta_{2},...,\zeta_{n}\right)
=\int_{W}\Omega\wedge\left(  \zeta_{1}\wedge\zeta_{2}\wedge...\wedge\zeta
_{n}\right)  \lrcorner\Omega.
\end{array}
\]
There is no instanton correction for the B-side. In particular identification
of these two Yukawa couplings would determine all genus zero Gromov-Witten
invariants of $M$. This is in fact the most well-known part of the mirror
symmetry conjecture.

\bigskip

\textbf{Batyrev's construction and the mirror theorem}

Before anyone can prove this COGP mirror conjecture, we need to know how to
construct $W$ from a given $M$. When $M$ is a Calabi-Yau hypersurface in a
Fano toric variety $X_{\Delta}$ as discussed on page \pageref{Toric Eg},
Batyrev \cite{Ba} proposes that $W$ can be constructed again as a hypersurface
in another Fano toric variety $X_{\nabla}$, possibly with a crepant
desingularization. In fact the two polytopes $\Delta$ and $\nabla$ are polar
dual to each other (see \cite{LV} for a physical explanation).

Batyrev verifies that the Hodge numbers for these pairs do satisfy the COGP
mirror conjecture, namely $\dim H^{p,q}\left(  M\right)  =\dim H^{n-p,q}%
\left(  W\right)  .$ To verify the rest of the COGP conjecture amounts to
computing the genus zero Gromov-Witten invariants of $M$. This problem can be
translated into computing certain coupled Gromov-Witten invariants of
$X_{\Delta}$. Since the toric variety $X_{\Delta}$ has a $T^{n}$ torus action,
the moduli space of stable curves in $X_{\Delta}$ also inherit such a $T^{n}$
action. In the foundational paper \cite{K2}, Kontsevich developes the
techniques needed to apply Bott's localization method to compute Gromov-Witten
invariant using this torus action. In \cite{Gi} Givental gives an argument to
compute these invariants using localization. A complete and detailed proof of
the COGP conjecture is given by Lian, Liu and Yau in \cite{LLY1} and
\cite{LLY2}. This result is also proven by Bertram by a different method later.

\bigskip

\textbf{Bershadsky-Cecotti-Ooguri-Vafa conjecture}

In \cite{BCOV1} and \cite{BCOV2}, the authors study higher genus Gromov-Witten
invariants in a Calabi-Yau threefold. They conjecture that their generating
functions satisfy a differential equation, so-called $tt^{\ast}$-equation. In
particular in the quintic threefold case, this is enough to predict all
Gromov-Witten invariants. The mathematical structure on the B-side is still
not completely understood yet. For this amazing conjecture, the reader is
referred to the original papers.

\bigskip

\textbf{Gopakumar-Vafa conjecture}

As we mentioned before, Gromov-Witten invariants for Calabi-Yau threefolds are
very difficult to compute. One of the reasons is that it is very difficult to
describe the moduli space of stable maps explicitly. For example the moduli
space of stable maps can be very complicated even though their image is a
single curve inside $M$. Moreover these invariants may not be integer and also
there are infinite number of such invariants in each homology class of $M$,
namely one for each genus.

Recently in \cite{GV1}, \cite{GV2}, Gopakumar and Vafa propose to study a
rather different space whose cohomology admits an $sl\left(  2\right)  \times
sl\left(  2\right)  $ action. Its multiplicities are called BPS numbers and
they should determine all Gromov-Witten invariants. The origin of this
conjecture comes from M-theory duality considerations in physics. First the
space they consider is just the space of holomorphic curves, instead of stable
maps, in $M$ together with flat $U\left(  1\right)  $ bundles over them. This
is just the moduli space of B-cycles $_{B}\mathcal{M}\left(  M\right)  $ on
$M$. Second, unlike Gromov-Witten invariants, these BPS numbers are always
integers! Third, there are only a finite number of BPS numbers for a given
homology class in $M$ because they are bounded by the dimension of the moduli
space. In particular Gromov-Witten invariants for various genus are not
independent! In a way the relationship between the BPS numbers and
Gromov-Witten invariant for Calabi-Yau threefolds is like the relationship
between Seiberg-Witten invariants and Donaldson invariants for four manifolds.
Unlike the Seiberg-Witten invariants, the BPS numbers are still very difficult
to compute. In fact they are not even well-defined at this moment.

The first issue is to compactify the moduli space $_{B}\mathcal{M}\left(
M\right)  $. When the curve $C$ in $M$ is smooth, then the space of flat
$U\left(  1\right)  $ bundles on $C$ is its Jacobian $J\left(  C\right)  $.
But when the curve becomes singular, especially nonreduced, there is no good
notion of compactified Jacobian. Presumably it will also depend on how $C$
sits inside $M$. The compactified moduli space $\overline{_{B}\mathcal{M}%
\left(  M\right)  }$, if it exists, would admit an Abelian variety fibration
obtained by forgetting the bundles. Roughly speaking the conjectured
$sl\left(  2\right)  \times sl\left(  2\right)  $ action on the cohomology of
$\overline{_{B}\mathcal{M}\left(  M\right)  }$ would correspond to the two
$sl\left(  2\right)  $ actions from the hard Lefschetz actions on the base and
the fiber of the Abelian variety fibration. It is still not clear how it works.

We now explain how to determine the Gromov-Witten invariants of curves in a
class $\beta\in H_{2}\left(  M,\mathbb{Z}\right)  $ in terms of the
$\mathbf{sl}\left(  2\right)  \mathbf{\times sl}\left(  2\right)  $ action
(\cite{GV2}). Let us denote the standard two dimensional representation of
$\mathbf{sl}\left(  2\right)  $ by $V_{1/2}$ and its $k^{th}$ symmetric power
by $V_{k/2}=S^{k}V_{1/2}$. Similarly the $\mathbf{sl}\left(  2\right)
\mathbf{\times sl}\left(  2\right)  $ representation $V_{j}\otimes V_{k}$ is
$V_{j,k}$ in our notation. We look at the moduli space of B-cycles $\left(
C,E\right)  $ with $\left[  C\right]  =\beta$ and $E$ a flat $U\left(
1\right)  $ bundle over $C$. Suppose that the conjectural $\mathbf{sl}\left(
2\right)  \mathbf{\times sl}\left(  2\right)  $ action exists on a suitable
cohomology theory of $_{B}\mathcal{M}\left(  M\right)  $. We decompose this
action as
\[
H^{\ast}\left(  _{B}\mathcal{M}\left(  M\right)  ,\mathbb{C}\right)  =\left[
V_{1/2,0}+2V_{0,0}\right]  \otimes\sum_{j,k}N_{j,k}^{\beta}V_{j,k}.
\]
We define integers $n_{\beta}^{g}\left(  M\right)  $ by the following,
\[
\sum_{g\geq0}n_{\beta}^{g}\left(  M\right)  \left[  V_{1/2}+2V_{0}\right]
^{\otimes g}=\sum_{j,k}N_{j,k}^{\beta}\left(  -1\right)  ^{2k}\left(
2k+1\right)  V_{j}\text{.}%
\]
These numbers are called the BPS numbers which count the numbers of BPS states.

The conjectural Gopakumar-Vafa formula expresses the Gromov-Witten invariants
of various genera curves in class $\beta$ in terms of these integers
$n_{\beta}^{g}\left(  M\right)  $'s. Their formula is%

\[
\sum_{\beta\neq0}\sum_{r\geq0}N_{\beta}^{r}\left(  M\right)  t^{2r-2}q^{\beta
}=\sum_{\beta\neq0}\sum_{r\geq0}n_{\beta}^{r}\left(  M\right)  \sum_{k>0}%
\frac{1}{k}\left(  2\sin\frac{kt}{2}\right)  ^{2r-2}q^{k\beta}\text{.}%
\]

For $\beta$ represented by a unique superrigid rational curve, namely a smooth
rational curve with normal bundle $O\left(  -1\right)  \oplus O\left(
-1\right)  $, this formula reduces to the multiple cover formula $N_{d\beta
}^{0}\left(  M\right)  =1/d^{3}$ in the genus zero case. For a superrigid
rational curve with one node, this formula is verified in \cite{BKL}. The
contribution of a super-rigid genus $g$ curve $C$ to $n_{d\left[  C\right]
}^{r}\left(  M\right)  $ has been computed through a certain range $r<R\left(
d,g\right)  $ by Bryan and Pandharipande in \cite{BP}. They verify the
integrality of these numbers which involves a non-trivial combinatorial
analysis of the number of degree $d$, etale covers of a smooth curve.

\bigskip

\textbf{Kontsevich homological mirror conjecture}

In \cite{K1} Kontsevich proposed a homological mirror symmetry conjecture: If
$M$ and $W$ are mirror manifolds, then a derived category of Lagrangians with
flat bundles in $M$ is isomorphic to the derived category of coherent sheaves
on $W$. These derived categories contain a lot of informations of $M$ and $W$.
Unfortunately we will not discuss much of those here because our emphasis is
on the geometric side; readers should consult his original paper \cite{K1} for
more details. This conjecture is verified by Polishchuk and Zaslow in
\cite{PZ} in dimension one, namely when $M$ and $W$ are elliptic curves. When
the dimension is bigger than one, the derived category of Lagrangians is not
even completely defined yet. One issue is about instanton corrections which
are needed in the definition of product structure in this category. They
involve a suitable count of holomorphic disks whose boundaries lie on
Lagrangians (see \cite{FOOO}). The COGP conjecture should follow from the
homological conjecture by considering diagonals in the products $M\times
\bar{M}$ and $W\times W$ as argued in \cite{K1}. It is because a holomorphic
disk in $M\times\bar{M}$ with boundary in the diagonal corresponds to two
holomorphic disks in $M$ with the same boundary, thus producing genus zero
holomorphic curves. This was explained to the author by Zaslow.

Another mirror conjecture of Kontsevich concerns symplectomorphism of $M$ and
holomorphic automorphism of $W$. As he suggested, one should enlarge the
collection of automorphisms of $W$ to automorphisms of the derived category of
coherent sheaves. Such automorphisms usually come from monodromy on the moduli
space of symplectic structures on $M$ and complex structures on $W$. And they
should have a mirror correspondence them. This is analyzed by Seidel and
Thomas \cite{ST} and Horja \cite{Ho} in some cases. We will explain how to
transform certain symplectic diffeomorphisms of $M$ to holomorphic
diffeomorphisms of $W$, and vice versa, in the semi-flat case in \cite{L3}.

Recently Kontsevich and Soibelman \cite{KS} have important results. They argue
that, under suitable assumptions, both the Fukaya category for $M$ and the
derived category of coherent sheaves on $W$ reduce to one category, which
should be realized as a Morse category for the base of a torus fibration in
the large complex structure limit.

\bigskip

\textbf{Strominger-Yau-Zaslow conjecture}

Next we discuss various geometric aspects of mirror symmetry. This starts with
the paper \cite{SYZ} by Strominger, Yau and Zaslow. They argue physically that
supersymmetric cycles on $M$ and $W$ should correspond to each other. The
manifold $W$ itself is the moduli space of points in $W$! The mirror of each
point, as a B-cycle in $W$, would be certain special Lagrangian with a flat
$U\left(  1\right)  $ connection, as a A-cycle in $M$. Because $\dim W=n$,
these special Lagrangian move in an $n$ dimensional family and therefore their
first Betti number equals $n$ (see section \ref{GeomCycle}). A natural choice
would be an $n$ dimension torus which moves and forms a special Lagrangian
torus fibration on $M$.

Consider at all those points in $W$ whose mirrors are the same special
Lagrangian torus in $M$ but with different flat $U\left(  1\right)  $
connections over it. Then this is just the dual torus, because the moduli
space of flat $U\left(  1\right)  $ connections on a torus is its dual. We
conclude that if $M$ and $W$ are mirror manifolds, then $W$ should have
special Lagrangian torus fibration and $M$ is the dual torus fibration.

It has been conjectured by various people including Kontsevich, Soibelman
\cite{KS}, Gross and Wilson, Todorov and Yau, that as $M$ and $W$ approach
certain large complex structure limit point, then they converge in the
Gromov-Hausdorff sense to affine manifolds. In the surface case, this is
verified by Gross and Wilson (see \cite{GW}). As we will see in the semi-flat
case, that the curvature is bounded in this limiting process and the affine
manifolds are the same via a Legendre transformation. We expect that the same
should hold true for general mirror pairs $M$ and $W,$ at least away from
singular fibers..

Fix $M$, a Calabi-Yau hypersurface in a Fano toric variety $X_{\Delta}$. As it
approaches the large complex structure limit point, i.e. the union of toric
divisors, then this Lagrangian torus fibration should be the one given by the
toric structure on $X_{\Delta}$. In particular this is a semi-flat situation.
The dual torus fibration is therefore the one given by $X_{\nabla}$ with the
polytope $\nabla$ the polar dual to $\Delta$. This offers an explanation to
why Batyrev construction should product mirror pairs. This (non-special)
Lagrangian torus fibration picture for Calabi-Yau inside toric varieties has
been verified by Gross \cite{Gr3}\footnote{The torus fibration constructed by
Gross is not Lagrangian.} and Ruan in \cite{R1},\cite{R3} (see also \cite{LV}
and \cite{Z}).

In order to understand the geometry of this moduli space $W$ of special
Lagrangians in $M$, for example to identify the L$^{2}$ metric with instanton
corrections on this moduli space with the Calabi-Yau metric on $W$, we need to
know how to count holomorphic disks. This point is still poorly understood.

\bigskip

\textbf{Vafa conjecture}

In \cite{V2} Vafa extended the SYZ analysis to general supersymmetric A-cycles
and B-cycles. Among other things, he wrote down a conjectural formula that
identifies the two correlation functions with instanton correction terms
included. This conjecture is verified in the special case when the B-cycle is
a line bundle on a semi-flat manifold in \cite{LYZ}.

However it is still not known how to correct the L$^{2}$ metric on the moduli
space to a Ricci flat metric, even conjecturally.

Roughly speaking, the mirror transformation between A-cycles and B-cycles
should be given by fiberwise Fourier transformation on the special Lagrangian
torus fibration, similar to the spectral cover construction in algebraic geometry.

For example if $M$ and $W$ are mirror manifolds with dual special Lagrangian
torus fibrations,
\begin{align*}
\pi &  :M\rightarrow D,\\
\pi &  :W\rightarrow D.
\end{align*}
Obviously $M\times M$ and $W\times W$ are also mirrors to each other. Now the
diagonal $\Delta_{W}$ in $W\times W$ is a complex submanifold and thus
determines a B-cycle in $W\times W$.

In Kontsevich picture, $\Delta_{W}$ corresponds to the identity functor, as a
kernel of the Fourier-Mukai transform, on the derived category of sheaves on
$W$. Therefore its mirror should correspond to the identity functor on the
Fukaya category of Lagrangians in $M$. The kernel for this identity function
is the diagonal $\Delta_{M}$ in $M\times\overline{M}$. It was explained to the
author by Thomas.

What is the mirror to $\Delta_{W}$ in SYZ picture?

We have $\pi_{W\times W}\left(  \Delta_{W}\right)  =\Delta_{B}\subset B\times
B $, the diagonal in $B$. Let $T_{b}=\pi_{W}^{-1}\left(  b\right)  $, then
\[
\Delta_{W}\cap\left(  T_{b}\times T_{b}\right)  =\Delta_{T_{b}}\subset
T_{b}\times T_{b}=\pi_{W\times W}^{-1}\left(  b,b\right)  \text{.}%
\]
The Fourier transformation on this fiber $T_{b}\times T_{b}$ would be those
flat $U\left(  1\right)  $ connections on $T_{b}\times T_{b}$ which are
trivial along $\Delta_{T_{b}}$ and this is the off-diagonal
\[
\widetilde{\Delta}_{T_{b}^{\ast}}=\left\{  \left(  y,-y\right)  \right\}
\subset T_{b}^{\ast}\times T_{b}^{\ast}.
\]
So the mirror to $\Delta_{W}$ in $M\times M$ is the union of $\widetilde
{\Delta}_{T_{b}^{\ast}}$'s over the diagonal $\Delta_{B}$ in $B\times B$.

If we map each point $y$ in the second factor of $T_{b}^{\ast}\times
T_{b}^{\ast}$ to its inverse $-y$, then the symplectic structure on $M\times
M$ would becomes the one in $M\times\overline{M}$. This can be verified
rigorously in the semi-flat case. Moreover the mirror to the B-cycle
$\Delta_{W}$ in $W\times W$ will be the A-cycle $\Delta_{M}$, the diagonal in
$M\times\overline{M}$.

\bigskip

\bigskip

\textbf{Donaldson, Uhlenbeck-Yau theorem and Thomas-Yau conjecture}

The relationship between the Kontsevich homological mirror symmetry conjecture
and SYZ geometric mirror symmetry conjecture is similar to the relationship
between the deRham-Singular cohomology and the Hodge cohomology of harmonic
forms for a manifold. To link the two, we need an analog of the Hodge theory
which represents each cohomology class by a unique harmonic form.

On the B-side, we need to find a canonical Hermitian metric on holomorphic
bundles over a complex submanifold in $W$, namely a deformed
Hermitian-Yang-Mills metric. For the usual Hermitian-Yang-Mills equation, this
is the result of Donaldson, Uhlenbeck and Yau which says that such equation
has a unique solution on any Mumford stable bundle. For the deformed
Hermitian-Yang-Mills equation near the large radius limit, i.e. $t\omega$ for
large $t$, the author shows in \cite{L1} that they can be solved on
irreducible Gieseker type stable bundles.

On the A-side, in \cite{Th2} Thomas formulates a stability notion for
Lagrangian submanifolds and conjectures that it should be related to the
existence of special Lagrangian in its Hamiltonian deformation class. Some
special cases are verified by Thomas and Yau \cite{TY} using mean curvature flow.

\bigskip

\textbf{Conjecture on} $sl\left(  2\right)  \times sl\left(  2\right)  $
\textbf{action on cohomology}

The hard Lefschetz theorem in K\"{a}hler geometry says that there is an
$\mathbf{sl}\left(  2\right)  $ action on the cohomology group of any compact
K\"{a}hler manifold. A generator for this $\mathbf{sl}\left(  2\right)  $
action is the hyperplane section on harmonic forms. Hubsch and Yau are the
first to search for the mirror of this $\mathbf{sl}\left(  2\right)  $ action
in \cite{HY}.

In the semi-flat case, there is another $\mathbf{sl}\left(  2\right)  $ action
whose generator is the variation of Hodge structure along the deformation
direction which shrinks the torus fiber (see \cite{L3}). Moreover these two
$\mathbf{sl}\left(  2\right)  $ actions commute with each other and thus
determine an $\mathbf{sl}\left(  2\right)  \times\mathbf{sl}\left(  2\right)
$ action, or equivalently an $\mathbf{so}\left(  4\right)  $ action.
Furthermore these two kinds of $\mathbf{sl}\left(  2\right)  $ actions get
interchanged under the mirror transformation. When the semi-flat Calabi-Yau
manifold is hyperk\"{a}hler, this $\mathbf{so}\left(  4\right)  $ action
imbeds into the natural hyperk\"{a}hler $\mathbf{so}\left(  5\right)  $ action
on its cohomology group \cite{L3}.

We conjecture: (1) if $M$ and $W$ are mirror manifolds, then their cohomology
groups admit an $\mathbf{so}\left(  4\right)  =\mathbf{sl}\left(  2\right)
\times\mathbf{sl}\left(  2\right)  $ action, where one of them is given by the
(possibly deformed) hard Lefschetz theorem, and the other one is related to
the variation of Hodge structure for a complex deformation which shrinks the
fibers of the special Lagrangian fibration. (2) Furthermore mirror
transformation flips these two kinds of $\mathbf{sl}\left(  2\right)  $
action, possibly after instanton corrections. (3) When $M$ is also a
hyperk\"{a}hler manifold, then this $\mathbf{so}\left(  4\right)  $ action
embeds inside the natural hyperk\"{a}hler $\mathbf{so}\left(  5\right)  $
action on its cohomology groups.

More generally suitable cohomology groups of the compactified moduli space of
supersymmetric A-cycles (or B-cycles) $\left(  C,E\right)  $ with $rank\left(
E\right)  =1$ admits an $\mathbf{so}\left(  4\right)  =\mathbf{sl}\left(
2\right)  \times\mathbf{sl}\left(  2\right)  $ action. The mirror
transformation between these moduli spaces for $M$ and $W$ flips the two
$sl\left(  2\right)  $ factors. When these cycles are just points or special
Lagrangian tori, this reduces to the previous part of the conjecture.

When the moduli space of B-cycles consists of curves, then this conjecture
overlaps with the first part of the Gopakumar-Vafa conjecture. It would be
interesting to know if there is an analog of the full Gopakumar-Vafa
conjecture in the general case.

\bigskip

\bigskip

\textbf{How to prove the mirror conjectures?}

In the ground breaking paper of Strominger, Yau and Zaslow [SYZ], a geometric
approach to proving the mirror conjecture is proposed. Namely the mirror
manifold $W$ of $M$ should be the moduli space of certain supersymmetric
cycles on $M$. Let me briefly explain their proposal here.

From string theory, the duality between $M$ and $W$ should have no instanton
correction at the large complex and K\"{a}hler structure limit. For example,
when $M$ is a Calabi-Yau hypersurface in projective space $\mathbb{P}^{n+1}$,
the large complex structure limit point $M_{0}$ is just the union of $n+2$
hyperplanes inside $\mathbb{P}^{n+1}$. Each connected component of the smooth
part of $M_{0}$ is a copy of $\left(  \mathbb{C}^{\times}\right)  ^{n}$ and it
has a semi-flat Calabi-Yau metric.

The most important step is to understand what the Calabi-Yau metric on $M$
looks like when $M$ is close to such a large complex structure limit point. It
is expected that on a large part on $M$, the Calabi-Yau metric on $M$
approximates the above semi-flat metric and therefore admits a special
Lagrangian fibration. When these special Lagrangian tori approach
singularities of $M_{0}$, it itself might develop singularities. In fact the
Calabi-Yau metric on $M$ should be close to some model Calabi-Yau metric
around singular special Lagrangian tori. Unfortunately the only such model
metric is in dimension two, namely the Ooguri-Vafa metric \cite{OV}. In
general degeneration of Calabi-Yau metrics is a difficult analytic problem.
However, in dimension two, Gross and Wilson \cite{GW} verified such a
prediction with the help of a hyperk\"{a}hler trick.

The second step would be to explain the mirror conjecture in the case of
semi-flat Calabi-Yau manifolds. This is partially done in \cite{LYZ} and
\cite{L3}. In fact understanding this semi-flat case enable us to enrich the
mirror conjecture as we saw in the last section. Notice that most of the
predictions in the semi-flat case do not involve instanton corrections. Then
we need to show that those duality properties in the semi-flat Calabi-Yau case
continue to hold true when the Calabi-Yau metric is only close to a semi-flat metric.

The third step is to write down many model Calabi-Yau metrics for
neighborhoods of singular special Lagrangian tori. By explicit calculations,
we hope to show that mirror conjecture holds true for such models. The last
step would be to glue the results from the previous two steps together. This
would then imply the mirror conjecture for Calabi-Yau manifolds near a large
complex structure limit point.

In conclusion, existence of certain types of special Lagrangian fibration is
responsible for the duality behavior in mirror symmetry. For rigid Calabi-Yau
manifolds, we know that traditional mirror conjecture cannot hold true for a
simple reason. There are some proposals of using higher dimensional Fano
manifolds as the mirror manifold. We will propose a link between these and the
SYZ conjecture in section \ref{RigidCY}.

We should mention that there are other approaches to prove the homological
mirror symmetry conjecture by Kontsevich and Soibelman \cite{KS}.

\bigskip\newpage

\bigskip

\section{\label{MSsemiflat}Mirror symmetry for semi-flat CY\ manifolds}

\textbf{Dimension reduction and Cheng-Yau manifolds}

As we mentioned before, mirror symmetry concerns any Calabi-Yau manifold which
admits a deformation into a large complex structure limit point. Near any such
point, the Calabi-Yau metric is expected to approximate a semi-flat Calabi-Yau
metric (see \cite{KS} for an approach to prove the Kontsevich homological
mirror symmetry conjecture in this case). In the semi-flat case, we can study
translation invariant solution to the complex Monge-Amp\`{e}re equation and
reduce the problem to a real Monge-Amp\`{e}re equation. Semi-flat Calabi-Yau
manifolds are introduced into mirror symmetry in the foundational paper
\cite{SYZ} and then further studied in \cite{H1}, \cite{Gr2} and \cite{LYZ}.
We start with a semi-flat Calabi-Yau manifold
\[
M=D\times i\mathbb{R}^{n}\subset\mathbb{C}^{n},
\]
where $D$ is a convex domain in $\mathbb{R}^{n}$. More generally we can
replace $D$ by an affine manifold and $M$ its tangent bundle. We only consider
translation invariant quantity, for example, a translation invariant
K\"{a}hler potential $\phi$ of the K\"{a}hler form $\omega=i\partial
\bar{\partial}\phi$ determines a Calabi-Yau metric on $M$ if and only if
$\phi$ is an elliptic solution of the real Monge-Amp\`{e}re equation
\[
\det\left(  \frac{\partial^{2}\phi}{\partial x^{j}\partial x^{k}}\right)
=const.
\]
The deepest result is the theorem of Cheng and Yau in \cite{CY} on the
Dirichlet problem for this equation: Given any bounded convex domain $D$ in
$\mathbb{R}^{n}$, there is a unique convex function $\phi$ defined on $D$
satisfying the real Monge-Amp\`{e}re equation equation and $\phi|_{\partial
D}=0$. We will call such manifold $D$ a \textit{Cheng-Yau manifold}. As we
will see below, the Legendre transformation of $\phi$ produces another
Cheng-Yau manifold $D^{\ast}$. This duality among Cheng-Yau manifolds is an
important part of the mirror symmetry at the large complex structure limit.

\bigskip

\textbf{The Fourier transformation and the Legendre transformation}

We can compactify imaginary directions by taking a quotient of $i\mathbb{R}%
^{n}$ by a lattice $i\Lambda$. That is we replace the original $M$ by
$M=D\times iT$ where $T$ is the torus $\mathbb{R}^{n}/\Lambda$ and the above
K\"{a}hler structure $\omega$ descends to $D\times iT$ and $M$ has the
following structures \cite{L3} (see also \cite{Gr2}):
\[%
\begin{array}
[c]{lll}%
\text{Structures on }M\text{:} &  & \\
\smallskip\text{ Riemannian metric} &  & g_{M}=\Sigma\phi_{jk}\left(
dx^{j}\otimes dx^{k}+dy^{j}\otimes dy^{k}\right) \\
\smallskip\text{ Complex structure} &  & \Omega_{M}=\Pi\left(  dx^{j}%
+idy^{j}\right)  =dz^{1}\wedge...\wedge dz^{n}\\
\smallskip\text{ Symplectic structure} &  & \omega_{M}=\Sigma\phi_{jk}%
dx^{j}\wedge dy^{k}.
\end{array}
\]

The projection map $\pi:M\rightarrow D$ is a special Lagrangian fibration on
$M$. The moduli space of special Lagrangian tori with flat $U\left(  1\right)
$ connections is just the total space of the dual torus fibration
$\pi:W\rightarrow D$. Moreover the $L^{2}$ metric on the moduli space is the
same as putting the dual metric on each dual torus fiber (see for example
\cite{LYZ}). This is what we call a Fourier transformation:
\[
dx^{j}\rightarrow dx^{j}\text{ and }dy^{j}\rightarrow dy_{j}.
\]
Here $y_{j}$'s are the dual coordinates on the dual torus fiber. Since $M$ is
the quotient of $TD$ by a lattice $\Lambda$, $W$ is naturally a quotient of
$T^{\ast}D$ by the dual lattice $\Lambda$. In particular, it has a canonical
symplectic form. Together with the metric, we obtain an almost complex
structure on $W$. Then $W$ has the following structures
\[%
\begin{array}
[c]{lll}%
\text{Structures on }W\text{:} &  & \\
\smallskip\text{ Riemannian metric} &  & g_{W}=\Sigma\phi_{jk}dx^{j}\otimes
dx^{k}+\phi^{jk}dy_{j}\otimes dy_{k}\\
\smallskip\text{ Complex structure} &  & \Omega_{W}=\Pi\left(  \phi_{jk}%
dx^{k}+idy_{j}\right)  =dz_{1}\wedge...\wedge dz_{n}\\
\smallskip\text{ Symplectic structure} &  & \omega_{W}=\Sigma dx^{j}\wedge
dy_{j}.
\end{array}
\]
Here $\left(  \phi^{jk}\right)  =\left(  \phi_{jk}\right)  ^{-1}$. Next we
perform a Legendre transformation along the base $D$, namely we consider a
change of coordinates $x_{k}=x_{k}\left(  x^{j}\right)  $ given by
\[
\frac{\partial x_{k}}{\partial x^{j}}=\phi_{jk},
\]
thanks to the convexity of $\phi$. Combining both transformations we have%

\[%
\begin{array}
[c]{llll}%
\text{Legendre transformation} & dx^{j} & \rightarrow &  dx_{j}=\Sigma
\phi_{jk}dx^{k}\\
\text{Fourier transformation} & dy^{j} & \rightarrow &  dy_{j}\text{.}%
\end{array}
\]
and%

\[%
\begin{array}
[c]{lll}%
\text{Structures on }W\text{:} &  & \\
\smallskip\text{ Riemannian metric} &  & g_{W}=\Sigma\phi^{jk}\left(
dx_{j}\otimes dx_{k}+dy_{j}\otimes dy_{k}\right) \\
\smallskip\text{ Complex structure} &  & \Omega_{W}=\Pi\left(  dx_{j}%
+idy_{j}\right)  =dz_{1}\wedge...\wedge dz_{n}\\
\smallskip\text{ Symplectic structure} &  & \omega_{W}=\Sigma\phi^{jk}%
dx_{j}\wedge dy_{k}.
\end{array}
\]

Moreover we have
\[
\phi^{jk}=\frac{\partial^{2}\psi}{\partial x_{j}\partial x_{k}}\text{,}%
\]
for some function $\psi\left(  x_{k}\right)  $, the Legendre transformation of
$\phi\left(  x^{j}\right)  $. It is obvious that the convexity of $\phi$ and
$\psi$ are equivalent to each other. Moreover
\[
\det\left(  \frac{\partial^{2}\phi}{\partial x^{j}\partial x^{k}}\right)
=C\text{ is equivalent to }\det\left(  \frac{\partial^{2}\psi}{\partial
x_{j}\partial x_{k}}\right)  =C^{-1}.
\]

Roughly speaking we have changed $T^{\ast}D$ to $TD^{\ast}$. Now the
similarity between $g_{M},\omega_{M}$ and $g_{W},\omega_{W}$ are obvious.
Therefore $W$ is also a semi-flat Calabi-Yau manifold. If we apply the same
procedure to $W$, we are going to recover $M$ itself, an inversion property.

If we use the complex polarization of our B-field as follows,
\begin{align*}
\omega^{n}  &  =i^{n}\Omega\bar{\Omega}\\
\operatorname{Im}e^{i\theta}\left(  \omega+i\beta\right)  ^{n}  &  =0\\
\operatorname{Im}e^{i\phi}\Omega &  =0\text{ on the zero section.}%
\end{align*}
Then under the mirror transformation to $W$, they become
\begin{align*}
\omega_{W}^{n}  &  =i^{n}\Omega_{W}\bar{\Omega}_{W}\\
\operatorname{Im}e^{i\theta}\Omega_{W}  &  =0\text{ on the zero section.}\\
\operatorname{Im}e^{i\phi}\left(  \omega_{W}+i\beta_{W}\right)  ^{n}  &  =0.
\end{align*}

We come back and look at only the Fourier transformation which does not depend
on the function $\phi$, but not the Legendre transformation which does depend
on $\phi$. On the tangent bundle $M=TD$, if we vary its symplectic structures
while keeping its natural complex structure fixed. This is equivalent to
having a family of solutions to the real Monge-Amp\`{e}re equation. Now on the
$W=T^{\ast}D$ side, the corresponding symplectic structure is unchanged,
namely $\omega_{W}=\Sigma dx^{j}\wedge dy_{j}$. But the complex structures on
$W$ vary because the complex coordinates on $W$ are given by $dz_{j}%
=\Sigma\phi_{jk}dx^{j}+idy_{j}$ which depends on particular solutions of the
real Monge-Amp\`{e}re equation.

To complete this identification, we need to complexified the symplectic
structures, $\omega_{M}^{\mathbb{C}}=\omega_{M}+i\beta_{M}$, by including a
B-field. We are using the real involution of $M$ to define a B-field on it.
The complexified real Monge-Amp\`{e}re equation is,
\[
\det\left(  \phi_{jk}+i\eta_{jk}\right)  =C,
\]
where $\beta_{M}=i\partial\bar{\partial}\eta\left(  x\right)  $ and $C$ is a
non-zero constant. If we write
\[
\theta_{jk}\left(  x\right)  =\phi_{jk}\left(  x\right)  +i\eta_{jk}\left(
x\right)  ,
\]
then the above equation becomes $\det\left(  \theta_{jk}\right)  =C$. The
above discussions can be extended to this case if we consider$\ $the following
modified transformation:%

\[%
\begin{array}
[c]{llll}%
\smallskip\text{ Legendre transformation} & dx^{j} & \rightarrow &
dx_{j}=\Sigma\phi_{jk}dx^{k}\\
\smallskip\text{ Fourier transformation} & dy^{j} & \rightarrow &
dy_{j}+\Sigma\eta_{jk}dx^{k}\text{.}%
\end{array}
\]
See \cite{L3} for details. Now we do have an identification between the moduli
space of complexified symplectic structures on $M$ with the moduli space of
complex structures on $W$ \footnote{Notice that it is important that the
B-fields satisfy the complexified Monge-Amper\'{e} equation instead of being a
harmonic two form.}. This identification is both a holomorphic map and an
isometry. Moreover the two Yukawa couplings are also identified under this
transformation. How about instanton correction? There are neither holomorphic
curves in $M$ or $W$, nor holomorphic disks whose boundaries lie on any fiber
or any section. In particular this verifies the COGP mirror conjecture in this case.

If we vary the complex structures on $M$ as follows: we define the new complex
structure on $M$ using the following holomorphic coordinates,
\[
z_{t}^{j}=x^{j}+ity^{j},
\]
for any $t\in\mathbb{R}_{>0}$. As $t$ goes to zero, this will approach a large
complex structure limit for $M$. The corresponding Calabi-Yau metric for any
$t$ is given by
\[
g_{t}=\Sigma\phi_{jk}\left(  \frac{1}{t}dx^{j}\otimes dx^{k}+tdy^{j}\otimes
dy^{k}\right)  .
\]
The given fibration $\pi:M\rightarrow D$ is always a special Lagrangian
fibration for each $t$. Moreover the volume form on $M$ is independent of $t$,
namely $dv_{M}=\omega^{n}/n!$. As $t$ goes to zero, the size of the fibers
shrinks to zero while the base gets infinitely large.

If we rescale the metric to $tg_{t}$, then the diameter of $M$ stays bound and
it converges in the Gromov-Hausdorff sense to the real $n$ dimension manifold
$D$ with the metric $g_{D}=\Sigma\phi_{jk}dx^{j}dx^{k}$ as $t$ approaches
zero. Moreover the sectional curvature stays bound in this limit. This is
consistent with the above mentioned expectations of what happens near the
large complex structure limit point.

In fact the variation of Hodge structures given by this family of variation of
complex structures on $M$ determines an $\mathbf{sl}\left(  2\right)  $ action
both on level of differential forms and on the level of cohomology groups of
$M$ (see \cite{L3}). Together with the hard Lefschetz $\mathbf{sl}\left(
2\right)  $ action, this produces an $\mathbf{sl}\left(  2\right)
\times\mathbf{sl}\left(  2\right)  $ action on $H^{\ast}\left(  M,\mathbb{C}%
\right)  $. Under the above mirror transformation, the VHS $\mathbf{sl}\left(
2\right)  $ action on cohomology groups of $M$ (resp. $W$) will become the
hard Lefschetz $\mathbf{sl}\left(  2\right)  $ action on cohomology groups of
$W$ (resp. $M$). This confirms our $\mathbf{sl}\left(  2\right)
\times\mathbf{sl}\left(  2\right)  $ conjecture in the semi-flat Calabi-Yau case.

In \cite{LYZ}, Yau, Zaslow and the author identifies the moduli space of
A-cycles on $M$ which are sections and B-cycles on $W$. We also identify their
correlation functions, thus verifies the Vafa conjecture in this case.

In conclusion, we can understand and prove many mirror symmetry phenomenons in
the semi-flat case via the modified Fourier transformation and the Legendre transformation.

\pagebreak 

\section{\label{RigidCY}Mirror symmetry for rigid CY manifolds}

Since the beginning of the mirror symmetry era, it is well-known that the
mirror symmetry conjecture cannot hold true for all Calabi-Yau manifolds. If
$M$ is a Calabi-Yau manifold with $H^{1}\left(  M,T_{M}\right)  =0$, then its
complex structure has no infinitesimal deformation. We call such $M$ a rigid
Calabi-Yau manifold. If $W$ is a mirror manifold of $M$, it would have
$H^{1}\left(  W,T_{W}^{\ast}\right)  =0$. In particular $W$ cannot be a
K\"{a}hler manifold, a contradiction. Such failure is caused by the
non-existence of a special Lagrangian fibration on $M$.

Nonetheless it is argued by Schimmrigk \cite{Sc}\cite{Schim2}, Candelas,
Derrick and Parkes \cite{CDP} and Vafa using conformal field theory that in
certain cases there should still be mirror manifolds which are certain higher
dimensional Fano manifolds. Batyrev and Borisov \cite{BB} give a systematic
way to construct candidates of such higher dimensional mirror manifolds using
toric geometry.

In this section we will explain the construction of these higher dimension
mirror manifolds and offer a partial answer to mirror symmetry phenomenons
using ideas from \cite{SYZ}. We begin with an example of rigid Calabi-Yau threefold.

\bigskip

\textbf{An example of rigid Calabi-Yau threefold}

Let $E_{0}$ be the Fermat cubic in $\mathbb{P}^{2}$:
\[
E_{0}=\left\{  \left[  z_{0},z_{1},z_{2}\right]  \in\mathbb{P}^{2}:z_{0}%
^{3}+z_{1}^{3}+z_{2}^{3}=0\right\}  =\mathbb{C/<}1,e^{\pi i/3}>.
\]
This is the only elliptic curve with a nontrivial $\mathbb{Z}_{3}$ action. If
$\omega$ is a nontrivial cubic root of unity, then the $\mathbb{Z}_{3}$ action
of $E_{0}$ is generated by
\[
\left[  z_{0},z_{1},z_{2}\right]  \rightarrow\left[  z_{0},\omega z_{1}%
,\omega^{2}z_{2}\right]  .
\]
This $\mathbb{Z}_{3}$ action has three fixed points: $\left[  1,0,0\right]  $,
$\left[  0,1,0\right]  $ and $\left[  0,0,1\right]  $. The diagonal action of
$\mathbb{Z}_{3}$ on $\left(  E_{0}\right)  ^{3}$ preserves its holomorphic
volume form and it has 27 fixed points. The quotient space
\[
M_{0}=\left(  E_{0}\right)  ^{3}/\mathbb{Z}_{3},
\]
is therefore a Calabi-Yau orbifold with 27 singular points, each model locally
on $\mathbb{C}^{3}/\mathbb{Z}_{3}$. Resolving $M_{0}$ by replacing each
singular point by a copy of $\mathbb{P}^{2}$ with normal bundle $O_{\mathbb{P}%
^{2}}\left(  -3\right)  $. We obtain a simply connected Calabi-Yau threefold
$M$ with $H^{1}\left(  M,T_{M}\right)  =0$, a rigid manifold.

It is easy to construct a special Lagrangian fibration on $\left(
E_{0}\right)  ^{3}$. We start with one on the elliptic curve $E_{0}$ which is
given by the projection to the imaginary part,
\[%
\begin{array}
[c]{cc}%
E_{0}= & \mathbb{C}/\mathbb{<}1,e^{\pi i/3}>\rightarrow\mathbb{R}\\
& z\rightarrow\operatorname{Im}z.
\end{array}
\]
This is a special Lagrangian fibration on $E_{0}$ and similarly we have one on
$\left(  E_{0}\right)  ^{3}$ using the product fibration. However such
fibration will not be invariant under the diagonal $\mathbb{Z}_{3}$ action and
therefore it cannot descend to one on $M_{0}$. This is consistent with the
non-existence of mirror for the rigid Calabi-Yau manifold $M$ from SYZ point
of view.

Physicists argued from the point of view in the conformal field theory that
there should still be a higher dimensional mirror manifold for $M$, say
$\bar{W}$. Concretely $\bar{W}$ is the quotient of the Fermat cubic
hypersurface in $\mathbb{P}^{8}$ by a $\mathbb{Z}_{3}$ action. For usual
mirror manifolds $M$ and $W$, their Hodge diamonds flip, in the sense that
$\dim H^{p,q}\left(  M\right)  =\dim H^{3-p,q}\left(  W\right)  $. In our
present situation, $\bar{W}$ has bigger dimension than $M$ and the flipped
Hodge diamond of $M$ sits in the middle of the Hodge diamond of $\bar{W}$,
i.e.
\[
\dim H^{p,q}\left(  M\right)  =\dim H^{3-p+2,q+2}\left(  \bar{W}\right)
\text{.}%
\]
Explicitly we have $\dim H^{1,1}\left(  M\right)  =36$, $\dim H^{2,1}\left(
M\right)  =0$ and $\dim H^{4,3}\left(  \bar{W}\right)  =36$ and $\dim
H^{3,3}\left(  \bar{W}\right)  =0$.

\bigskip

\textbf{A cubic fourfold as the mirror of a K3 surface}

To get a closer look at what happens, we consider the two dimensional case
first. The simplest example is the Kummer K3 surface. The map which sends $x$
to $-x$ generates a $\mathbb{Z}_{2}$ action on $\left(  E_{0}\right)  ^{2}$,
or any complex two torus. The quotient $\left(  E_{0}\right)  ^{2}%
/\mathbb{Z}_{2}$ has 16 singular points, each is an ordinary double point.
Namely it is locally modelled on $\mathbb{C}^{2}/\mathbb{Z}_{2}$ and it is
called an $A_{1}$ singularity. Blowing up will replace each singular point by
a copy of $\mathbb{P}^{1}$ with normal bundle $O_{\mathbb{P}^{1}}\left(
-2\right)  $ and we obtain a smooth simply connected Calabi-Yau surface called
the Kummer K3 surface.

To construct a special Lagrangian fibration on the Kummer K3 surface, we start
with the one on $E_{0}$ and therefore $\left(  E_{0}\right)  ^{2}$ as before.
A simple observation: this special Lagrangian fibration on $\left(
E_{0}\right)  ^{2}$ is invariant under the $\mathbb{Z}_{2}$ action and
therefore descends to one on $\left(  E_{0}\right)  ^{2}/\mathbb{Z}_{2}$.
Using twistor transformation to an elliptic fibration problem, one can show
that this special Lagrangian fibration structure continue to exist on the
Kummer K3 surface.

Next we consider $M_{0}=\left(  E_{0}\right)  ^{2}/\mathbb{Z}_{3}$ where the
first $\mathbb{Z}_{3}$ action is as before and the second one is its square.
Now $M_{0}$ is a Calabi-Yau orbifold with 9 isolated singular points of type
$A_{2}$. Resolving $M_{0}$ by blowing up these singularities will replace each
singular point by two smooth rational curves meeting a point. The normal
bundle of each rational curve is $O_{\mathbb{P}^{1}}\left(  -2\right)  . $
This again produces a K3 surface $M$.

The natural special Lagrangian fibration on $\left(  E_{0}\right)  ^{2}$ is
not invariant under the $\mathbb{Z}_{3}$ action and therefore does not induce
one on $M_{0}$ nor $M$. According to the recipe by Batyrev and Borisov, the
mirror of such $M$ should be the Fermat cubic fourfold in $\mathbb{P}^{5}$:
\[
\bar{W}=\left\{  z_{0}^{3}+z_{1}^{3}+z_{2}^{3}+z_{3}^{3}+z_{4}^{3}+z_{5}%
^{3}=0\right\}  \subset\mathbb{P}^{5}\text{.}%
\]

We need to come back to understand why the special Lagrangian fibration on
$\left(  E_{0}\right)  ^{2}$ does not work. After twistor rotation, each fiber
in $\left(  E_{0}\right)  ^{2}$ is a genus one Riemann surface but its image
in $\left(  E_{0}\right)  ^{2}/\mathbb{Z}_{3}$ is no longer embedded. In fact
it intersects itself transversely at one point. Namely it is a singular
Riemann surface of genus two! A generic deformation of this supersymmetric
cycle will smooth out the node. Such nearby cycles will no longer come from
the special Lagrangian fibration on $\left(  E_{0}\right)  ^{2}$.

The moduli space of these supersymmetric A-cycles is not going to produce us
the usual mirror manifold of $M$, which is usually $M$ itself because it is
hyperk\"{a}hler. Instead it is the compactified Jacobian $\mathcal{J}$ for
genus two curves in the K3 surface $M$. Birationally this hyperk\"{a}hler
manifold $\mathcal{J}$ is the Hilbert scheme of 2 points on $M$. In fact
$\mathcal{J}$ is also birational to the Fano variety of lines on the cubic
fourfold $\bar{W}$. In terms of our languages, this is a moduli space of
B-cycles $\left(  C,E\right)  $ on $\bar{W}$ such that $C$ is a degree one
curve in $\bar{W}$ and $E$ is a unitary flat line bundle over $C$ which is
necessarily trivial. Therefore the Batyrev-Borisov mirror $\bar{W}$ and the
SYZ mirror $\mathcal{J}$ of $M$ are related as one being the moduli space of
B-cycles on another one. It is curious to know how much of these continue to
happen for more general cases. The author thanks B. Hassett for many
discussions on these matters.

\bigskip

\textbf{Generalized Calabi-Yaus and Batyrev-Borisov construction}

We are going to describe the Batyrev-Borisov construction of higher
dimensional candidates of mirror manifolds to complete intersection Calabi-Yau
manifolds inside toric varieties. In \cite{Sc} Schimmrigk described a class of
Fano manifolds, called generalized Calabi-Yau manifolds, that certain part of
its conformal field theory behaves as if it comes from a Calabi-Yau manifold.
For example it has a mirror which is again another generalized Calabi-Yau
manifold. A systematic way to construct these mirrors of generalized
Calabi-Yau manifolds using toric geometry is given by Batyrev and Borisov in
\cite{BB}. Moreover, given any Calabi-Yau manifold $M$ of dimension $n$ which
is a complete intersection of $r$ hypersurfaces in a Fano toric variety, there
is a generalized Calabi-Yau manifold $\bar{M}$ of dimension $n+2\left(
r-1\right)  $ associated to it. When $M$ and $W$ are mirror manifolds in the
usual sense, their associated generalized Calabi-Yau manifolds $\bar{M}$ and
$\bar{W}$ will continue to be mirror manifolds in the generalized sense.

When the complete intersection Calabi-Yau manifold $M$ is rigid, its mirror
manifold $W$ does not exist. However we still have $\bar{W}$, a mirror to the
generalized Calabi-Yau manifold $\bar{M}$, which we regard as the higher
dimensional mirror to $M$.

Now we describe the class of generalized Calabi-Yau manifolds $\bar{M}$, they
are hypersurfaces in Fano toric varieties $X_{\Delta}$. Recall that the zero
set of a section of $K_{X_{\Delta}}^{-1}$ is a Calabi-Yau manifold, possibly
singular. Its mirror is another Calabi-Yau which is the zero set of a section
of $K_{X_{\triangledown}}^{-1}$, here the polytope $\nabla$ is the polar dual
to $\Delta$. If we can write
\[
K_{X_{\Delta}}^{-1}=O\left(  r\right)  =O\left(  1\right)  ^{\otimes r}%
\]
the $r^{th}$ tensor power of a line bundle $O\left(  1\right)  $, then the
zero set of a section of $O\left(  1\right)  $ is called a generalized
Calabi-Yau manifold $\bar{M}$, and we have
\[
K_{\bar{M}}^{-1}=O\left(  r-1\right)  \text{.}%
\]
In particular $\bar{M}$ is a Fano manifold. When this happens, the
anti-canonical line bundle $K_{X_{\triangledown}}^{-1}$ of $X_{\triangledown}
$ also have the same property, its corresponding section will determine
another generalized Calabi-Yau manifold $\bar{W}$, we call it the mirror of
$\bar{M}$. We define the integer $n$ as
\[
\dim\bar{M}=\dim\bar{W}=n+2\left(  r-1\right)  \text{.}%
\]
To study mirror symmetry for these generalized Calabi-Yau manifolds, we should
not include all complex submanifolds. For example points should not be
considered as B-cycles. How about A-cycles? We should probably look at all
Lagrangian submanifolds, but we do not have a good notion of a special
Lagrangian submanifold. For the homological mirror symmetry, we can still
discuss the Fukaya category of Lagrangians in $\bar{M}$, or $\bar{W}$. We
would look for a category of certain coherent sheaves on $\bar{M}$, or
$\bar{W}$ so that this category behaves as if it is a category of coherent
sheaves on a Calabi-Yau manifold of dimension $n$. For instance $Ext^{k}%
\left(  S_{1},S_{2}\right)  ^{\ast}=Ext^{n-k}\left(  S_{2},S_{1}\right)  $.
The natural question would be to compare these two types of categories.

\bigskip

\textbf{From Calabi-Yaus to generalized Calabi-Yaus}

We use the above construction of mirror for generalized Calabi-Yau manifolds
to obtain higher dimensional mirrors for rigid Calabi-Yau manifolds which are
complete intersections in Fano toric varieties.

We first recall that when $M$ is a Calabi-Yau hypersurface in a Fano toric
variety $X_{\Delta}$, it can be deformed to the most singular Calabi-Yau which
is the union of toric divisors in $X_{\Delta}$. This is the large complex
structure limit. It is expected that this is the origin of special Lagrangian
fibration on $M$. The same idea works for Calabi-Yau complete intersections in
$X_{\Delta}$ provided that there are enough deformations. This is the case if
\[
M=\bigcap_{i=1}^{r}D_{i}%
\]
with each $D_{i}$ semi-ample and $\Sigma D_{i}=-K_{X_{\Delta}}$.

If some of the $D_{i}$'s are not semi-ample then the Calabi-Yau manifold $M$
may not be able to deform to a large complex structure limit point. And we do
not expect to find a special Lagrangian fibration on $M$ coming from its
embedding inside $X_{\Delta}$. This is precisely what happens for the rigid
Calabi-Yau threefold $\left(  E_{0}\right)  ^{3}/\mathbb{Z}_{3}$. Of course
this will not happen if $r=1$ because $D=-K_{X_{\Delta}}$ is ample for Fano variety.

The collection of effective divisors $D_{i}$'s in $X_{\Delta}$ determines a
hypersurface $\bar{M}$ in the projective bundle $\mathbb{P}\left(
\oplus_{i=1}^{r}O\left(  D_{i}\right)  \right)  $ over $X_{\Delta}$.
\[
\bar{M}\subset\mathbb{P}\left(  \oplus_{i=1}^{r}O\left(  D_{i}\right)
\right)  .
\]
This is because a section of $O\left(  D_{i}\right)  $ over $X_{\Delta}$ can
be identified as a function on the total space of the line bundle $O\left(
-D_{i}\right)  $ which is linear along fibers. If we projectivize each fiber
of the direct sum of $O\left(  D_{i}\right)  $'s, then this function will
become a section of $O\left(  1\right)  ,$ the dual of the tautological line
bundle. This gives all sections of $O\left(  1\right)  $, namely
\[
\oplus_{i=1}^{r}H^{0}\left(  X_{\Delta},O\left(  D_{i}\right)  \right)  \cong
H^{0}\left(  \mathbb{P}\left(  \oplus_{i=1}^{r}O\left(  D_{i}\right)  \right)
,O\left(  1\right)  \right)  \text{.}%
\]
Note that $\bar{M}$ is a generalized Calabi-Yau manifold with $\dim\bar
{M}=\dim M+2r-2$ and $h^{p,q}\left(  M\right)  \leq h^{p+r-1,q+r-1}\left(
\bar{M}\right)  $. If we restrict the projective bundle morphism to $\bar{M}$,
we obtain a surjective morphism
\[
\pi:\bar{M}\rightarrow X_{\Delta}\text{.}%
\]
The fiber over a point in $M$ (resp. outside $M$) is a copy of $\mathbb{P}%
^{r-1}$ (resp. $\mathbb{P}^{r-2}$). Therefore $Y=\pi^{-1}\left(  M\right)
\subset\bar{M}$ is a $\mathbb{P}^{r-1}$-bundle over $M$,
\[%
\begin{array}
[c]{ccc}%
\mathbb{P}^{r-1} & \rightarrow &  Y\\
&  & \,\downarrow{\small \pi}_{Y}\\
&  & M.
\end{array}
\]
In fact $M$ is the moduli space of $\mathbb{P}^{r-1}$ in the fiber homology
class in $\bar{M}$.

\bigskip

\textbf{A construction of Lagrangian submanifolds in }$\bar{M}$

Since $\bar{M}$ has higher dimension than $M$, it is a very interesting
problem to compare the symplectic geometries between $M$ and $\bar{M}$. The
author benefits from discussions with Seidel.

Using the method of vanishing cycles, we give a new construction of a
Lagrangian submanifold in $\bar{M}^{n+2r-2}$ from one in $M^{n}$: Recall that
the hypersurface $\bar{M}\subset\mathbb{P}\left(  \oplus_{i=1}^{r}O\left(
D_{i}\right)  \right)  $ is defined by $s_{1}+\ldots+s_{r-1}+s_{r}=0$. Notice
each $s_{j}$ is only well-defined up to nonzero multiple. If we rescale
$s_{r}$ to zero then
\[
\bar{M}_{0}=\left\{  s_{1}+\ldots+s_{r-1}+s_{r}=0\right\}  \subset
\mathbb{P}\left(  \oplus_{i=1}^{r}O\left(  D_{i}\right)  \right)
\]
is a cone over $X_{\Delta}$ of
\[
\bar{N}=\left\{  s_{1}+\ldots+s_{r-1}=0\right\}  \subset\mathbb{P}\left(
\oplus_{i=1}^{r-1}O\left(  D_{i}\right)  \right)  .
\]
The fiber of $\bar{M}_{0}\rightarrow X_{\Delta}$ over a point in
$\bigcap_{i=1}^{r-1}D_{i}$ (resp. in $X_{\Delta}\backslash\bigcap_{i=1}%
^{r-1}D_{i}$) is a copy of $\mathbb{P}^{r-1}$ (resp. $\mathbb{P}^{r-2}$).
Moreover $\bar{M}_{0}$ is singular along those cone points over $\bigcap
_{i=1}^{r-1}D_{i}$. When $\bar{M}$ degenerates to $\bar{M}_{0}$ as we rescale
$s_{r}$ to zero, the vanishing cycles is a family of spheres $S^{2r-3}$
parametrized by $\bigcap_{i=1}^{r-1}D_{i}=Sing\left(  \bar{M}_{0}\right)  $.
However, over any point in $M\subset\bigcap_{i=1}^{r-1}D_{i}$, such a
$S^{2r-3}$ in fact shrink to a point! This is because the degeneration from
$\bar{M}$ to $\bar{M}_{0}$ is not semi-stable and the base locus is given by
$M$. Being in the base locus means that such a point does not move during the
degeneration process, thus we never create a vanishing cycle for it - the
vanishing of the vanishing cycle.

Suppose $L$ is a Lagrangian submanifold in $M=\bigcap_{i=1}^{r}D_{i}$. We
assume that there is a Lagrangian submanifold $\hat{L}$ in $\bigcap
_{i=1}^{r-1}D_{i}$ with $\partial\hat{L}=L$. Such a $\hat{L}$ can be
constructed, for example, when $L$ is a vanishing cycle for $M$ as $D_{r}$
deforms. The above family of $S^{2r-3}$ over $\hat{L}$ which shrink to points
over its boundary $L$. The total space $\bar{L}$ is a \textit{closed}
submanifold of half dimensional in $\bar{M}$. With more care, one should be
able to verify that $\bar{L}$ is a Lagrangian submanifold in $\bar{M}$.

A natural and important problem is to compare the Floer homology groups
$HF\left(  M;L,L^{\prime}\right)  $ and $HF\left(  \bar{M}^{\prime};\bar
{L},\bar{L}^{\prime}\right)  $.

\bigskip

\textbf{Complex geometries between }$M$ \textbf{and }$\bar{M}$

For any coherent sheaf $S$ on $M$ we can easily construct an another coherent
sheaf on $\bar{M}$ as $\bar{S}=\iota_{\ast}\circ\pi_{Y}^{\ast}\left(
S\right)  $ where $\iota:Y\subset\bar{M}$ is the natural inclusion. An
important problem is to recover the derived category of coherent sheaves on
$M$ from the one on $\bar{M}$ together with $O\left(  1\right)  =K_{\bar{M}%
}^{1/\left(  r-1\right)  }$.

Similarly we can construct a homomorphism between their Chow groups,
\begin{align*}
CH_{k}\left(  M\right)   &  \rightarrow CH_{k+r-1}\left(  Y\right)
\rightarrow CH_{k+r-1}\left(  \bar{M}\right) \\
C  &  \rightarrow\pi_{Y}^{-1}C.
\end{align*}
Here the first morphism is the pullback map and the second morphism simply
regards $\pi_{Y}^{-1}C$ as a algebraic cycle in $\bar{M}$ via the embedding
$Y\subset\bar{M}$. In terms of cohomology groups, we have
\[
H^{l}\left(  M\right)  \overset{\alpha}{\rightarrow}H^{l}\left(  Y\right)
\overset{\beta}{\underset{\cong}{\rightarrow}}H^{l+2r-2}\left(  \bar{M}%
,\bar{M}\backslash Y\right)  \overset{\gamma}{\rightarrow}H^{l+2r-2}\left(
\bar{M}\right)  .
\]
Here $\alpha=\pi_{Y}^{\ast}$ is the pullback morphism, $\beta$ is the Thom
isomorphism and $\gamma$ is natural homomorphism in the long exact sequence
for the pair $\left(  \bar{M},\bar{M}\backslash Y\right)  $. We expect that
these morphisms preserve Hodge structures, in particular, the image of
$H^{p,q}\left(  M\right)  $ lies inside $H^{p+r-1,q+r-1}\left(  \bar
{M}\right)  $.

It seems more natural to consider $H^{\ast,\ast}\left(  \bar{M},\bar
{M}\backslash Y\right)  $ than $H^{\ast,\ast}\left(  \bar{M}\right)  .$ For an
arbitrary generalized Calabi-Yau manifold $\bar{M}$ with $K_{\bar{M}}%
^{-1}=O\left(  r-1\right)  $, we define $M$ to be the moduli space of
$\mathbb{P}^{r-1}$ representing the class $c_{1}\left(  O\left(  1\right)
\right)  ^{r-1}$ and $Y$ is the image of the universal variety inside $\bar{M}
$. For example, in the case of cubic sevenfold (or cubic fourfold), we have
$Y=\bar{M}$ and therefore $H^{\ast,\ast}\left(  \bar{M},\bar{M}\backslash
Y\right)  =H^{\ast,\ast}\left(  \bar{M}\right)  $.

When $\bar{M}$ is associated to a Calabi-Yau manifold $M$ as above, then we
have the Thom isomorphism,
\[
H^{l+2r-2}\left(  \bar{M},\bar{M}\backslash Y\right)  \cong H^{l}\left(
Y\right)  .
\]
So we can simply look at $Y$, a $\mathbb{P}^{r-1}$-bundle over $M$, rather
than the pair $\left(  \bar{M},\bar{M}\backslash Y\right)  $. For example, if
$E$ is a holomorphic vector bundle over $M$, then $\pi_{Y}^{\ast}E$ is a
holomorphic vector bundle on $Y$ and we have an isomorphism,
\[
H^{k}\left(  Y,O\left(  \pi_{Y}^{\ast}E\right)  \right)  \cong H^{k}\left(
M,E\right)  ,
\]
using the Leray spectral sequence. We would like to extend $\pi_{Y}^{\ast}E$
over $Y$ to a holomorphic bundle $\bar{E}$ over the pair $\left(  \bar{M}%
,\bar{M}\backslash Y\right)  $ (up to quasi-isomorphism) and obtain an analog
of the Thom isomorphism twisted by $\bar{E}$. More generally one might want to
perform these comparisons on the derived categories of coherent sheaves on $M$
and $\left(  \bar{M},\bar{M}\backslash Y\right)  $ (see the remark below).

We come back to discuss the mirror of generalized Calabi-Yau manifolds. If
$D_{i}$'s are semi-ample so that $M$ has a mirror manifold $W$, then $W$ is
also a complete intersection in a toric variety $X_{\nabla}$ where $\nabla$ is
the polar dual of the polytope $\Delta$. So there is also a generalized
Calabi-Yau manifold $\bar{W}$ associated to $W$. It can be shown that $\bar
{W}$ coincides with the Batyrev-Borisov's construction of the mirror for the
generalized Calabi-Yau manifold $\bar{M}$.

When $D_{i}$'s are not semi-ample, the mirror of $M$ might not exist but such
higher dimensional mirror $\bar{W}$ always exists. The generalization of
mirror symmetry conjecture to this situation should be the comparison of the
symplectic geometry on $M$ to certain part of the complex geometry on $\bar
{W}$.

\bigskip

\textbf{Rigid Calabi-Yau manifolds as moduli spaces}

A lesson we learned from the paper \cite{SYZ} is that the mirror manifold of
$M$ can be realized as a compactified moduli space of certain supersymmetric
A-cycles on $M$. In this case they are special Lagrangian tori with flat
connections on them. Roughly speaking these cycles sweep over $M$ once. This
explains why the geometry of $M$ is reflected in the geometry of this moduli
space $W$. A similar and simpler example is $M$ itself is a moduli space of
supersymmetric B-cycles in $M$, namely the moduli space of points in $M$.

In general the complex (resp. symplectic) geometry of $M$ is closely related
to the complex (resp. symplectic) geometry of any moduli space of B-cycles in
$M$. And the symplectic (resp. complex) geometry of $M$ is closely related to
the symplectic (resp. complex) geometry of any moduli space of A-cycles in $M$.

In fact any Calabi-Yau manifold $M$ in $X_{\Delta},$ which is a complete
intersection of $r$ hypersurface, is itself a moduli space of supersymmetric
B-cycles in $\bar{M}$. Each cycle $\left(  C,E\right)  $ in this moduli space
consists of a projective space $\mathbb{P}^{r-1}=C$ and the trivial line
bundle $E$ over it $C$. Moreover each of these $C$'s is a fiber of the
projective bundle over $X_{\Delta}$. This suggests an explanation for why $M $
and $\bar{W}$ should be mirror manifolds. For example, for the particular K3
surface whose mirror is a cubic fourfold as before, we have $r=2$ and that is
why we looked at the moduli space of $\mathbb{P}^{1}$. If we look at the
moduli space of those $\mathbb{P}^{1}$'s which are of degree one, namely the
Fano variety of lines on the cubic fourfold then it is, at least birationally,
the moduli space of those special Lagrangians on the K3 surface which descend
from the special Lagrangian fibration of the four torus. Since this moduli
space is a hyperk\"{a}hler manifold, it is self-mirror. Hence this provides a
possible explanation of mirror symmetry for generalized Calabi-Yau manifolds
from the SYZ picture in this case.

\bigskip

\textbf{A remark on Yau's program on noncompact CY manifolds}

A generalized Calabi-Yau $\bar{M}$ might not look like a Calabi-Yau manifold,
however its complement $X_{\Delta}\backslash\bar{M}$ does. Namely $X_{\Delta
}\backslash\bar{M}$ admits a complete Ricci flat K\"{a}hler metric if $\bar
{M}$ is a Calabi-Yau hypersurface or a generalized Calabi-Yau K\"{a}%
hler-Einstein manifold. We expect that this complement $X_{\Delta}%
\backslash\bar{M}$ also admit a special Lagrangian fibration and mirror
symmetry. Yau first discusses the problem of constructing K\"{a}hler-Einstein
metrics in the noncompact setting in \cite{Yau ICM talk} (see also \cite{Yau
Nonlinear}).

To begin, we suppose that $M$ is a Calabi-Yau manifold of dimension $n$ which
is an anti-canonical divisor in a Fano toric variety $X_{\Delta}$. When $M$
approaches the large complex structure limit, which is the union of $n$
dimensional strata of $X_{\Delta}$, then we expect from SYZ proposal
\cite{SYZ} that $M$ should admit a special Lagrangian fibration which looks
like the toric fibration on these strata in $X_{\Delta}$. The dual fibration
would product the mirror manifold of $M$ (see also \cite{LV}).

Not just $M$ admits a Ricci-flat K\"{a}hler metric, its complement $X_{\Delta
}\backslash M$ also admits a complete Ricci-flat K\"{a}hler metric by the work
of Tian and Yau \cite{TY1} and also Bando and Kobayashi \cite{BK} as the first
step towards Yau's problem on noncompact CY manifolds \cite{Yau ICM talk}.
When $M$ approaches the large complex structure limit, then $X_{\Delta
}\backslash M$ would becomes $\left(  \mathbb{C}^{\times}\right)  ^{n+1}$.
This open manifold certainly admits a special Lagrangian fibration with
respect to the standard flat metric. As an extension to the SYZ conjecture, we
expect that $X_{\Delta}\backslash M$ would admit a special Lagrangian
fibration. Moreover this fibration is compatible with the one on $M$, namely
when the fiber torus $T^{n+1}$ in $X_{\Delta}\backslash M$ goes to infinity,
which is $M$, then a circle direction collapses and it becomes a special
Lagrangian torus $T^{n}$ for the fibration in $M$. The mirror conjectures for
$M\subset X_{\Delta}$ and $W\subset X_{\triangledown}$ should also be extended
to $X_{\Delta}\backslash M$ and $X_{\triangledown}\backslash W.$

Now we consider a generalized Calabi-Yau manifold $\bar{M}$, which is a Fano
hypersurface in $X_{\Delta}$ with
\[
K_{X_{\Delta}}^{-1}=O\left(  r\bar{M}\right)  \text{.}%
\]
In \cite{TY2} and also \cite{BaK}, the authors prove that $X_{\Delta
}\backslash\bar{M}$ admits a complete Ricci-flat K\"{a}hler metric provided
that $\bar{M}$ is K\"{a}hler-Einstein. It would be useful if we can construct
a special Lagrangian fibration on $X_{\Delta}\backslash\bar{M}$ and understand
its behaviors near infinity. We also expect that its dual fibration should
coincide with the special Lagrangian fibration on $X_{\triangledown}%
\backslash\bar{W}$. This is a very difficult analytic problem but its solution
would provide us knowledge about mirror symmetry for generalized Calabi-Yau
manifolds from SYZ point of view.

For example when $\bar{M}$ is a point, say the north pole, in $X_{\Delta
}=\mathbb{P}^{1}$, then $X_{\Delta}\backslash\bar{M}=\mathbb{C}$ has a special
Lagrangian fibration by lines with constant real part. Its dual fibration
coincides with the corresponding fibration on $X_{\triangledown}\backslash
W=\mathbb{C}$. An anti-canonical Calabi-Yau hypersurface $M$ in $\mathbb{P}%
^{1} $ has two points, say the north pole and one other point $p$. Its
complement $X_{\Delta}\backslash M$ is isomorphic to $\mathbb{C}^{\times}$ and
it has a special Lagrangian fibration by $S^{1}$. If we move $p$ toward the
north pole, then $M$ becomes $2\bar{M}$, a non-reduced scheme. Moreover the
limiting special Lagrangian fibration will be the one above on $\mathbb{C}$
given by straight lines.

In general we should treat $\bar{M}\subset X_{\Delta}$ to have multiplicity
$r$. This is because the zero sets of a family of sections of $K_{X_{\Delta}%
}^{-1}=O\left(  r\bar{M}\right)  $, which are Calabi-Yau manifolds, can
converge to a non-reduced scheme $r\bar{M}$. The scheme structure on $\bar{M}
$ does depend on the family.

\section{Appendix: Deformations of SUSY B-cycles}

We are going to study the deformation theory of B-cycles and structures of
their moduli spaces. If we forget the deformed Hermitian-Yang-Mills equations,
namely we consider deformations of a complex submanifold $C$ together with a
holomorphic bundle $E$ over it, then they are parametrized infinitesimally by
the cohomology group $Ext_{O_{M}}^{1}\left(  \iota_{\ast}E,\iota_{\ast
}E\right)  $. To understand the deformation theory when the deformed
Hermitian-Yang-Mills equation is included, we would need to have differential
form representatives for classes in $Ext_{O_{M}}^{1}\left(  \iota_{\ast
}E,\iota_{\ast}E\right)  $. We should see that
\[
Ext_{O_{M}}^{1}\left(  \iota_{\ast}E,\iota_{\ast}E\right)  =H^{1}\left(
C,\mathbb{E}\right)  \times_{H^{1}\left(  C,T_{C}\right)  }H^{0}\left(
C,N_{C/M}\right)  ,
\]
and we can represent elements on the right hand side by differential forms.
First we consider only deformations of the deformed Hermitian-Yang-Mills
equation while $C$ is kept fixed.

\bigskip

\textbf{Deformation of deformed Hermitian-Yang-Mills bundles}

To begin we consider infinitesimal deformations of B-cycles $\left(
C,E\right)  $ which fix the complex submanifold $C\subset M$. If we forget the
Hermitian metric on $E$ and only deform the holomorphic structure of the
vector bundle $E$, then they are parametrized by the cohomology group
$H^{1}\left(  C,End\left(  E\right)  \right)  $ infinitesimally. There are
nontrivial \cite{Th1} obstruction classes in $H^{2}\left(  C,End\left(
E\right)  \right)  $ for these infinitesimal deformations to come from a
honest deformation.

If we vary Hermitian-Yang-Mills connections on $E$, then it is well-known that
they are parametrized by harmonic $\left(  0,1\right)  $ forms with valued in
$ad\left(  E\right)  $. By Dolbeault-Hodge theorem this is isomorphic to
$H^{1}\left(  C,End\left(  E\right)  \right)  $. It is also not difficult to
see that infinitesimal variations of deformed Hermitian-Yang-Mills bundles on
$C$ are parametrized by $B\in\Omega^{0,1}\left(  C,ad\left(  E\right)
\right)  $ satisfying
\begin{align*}
\bar{\partial}B  &  =0,\\
\operatorname{Im}e^{i\theta}\left(  \omega^{\mathbb{C}}+F\right)  ^{m-1}%
\wedge\partial B  &  =0.
\end{align*}
In general we have the following definition (\cite{LYZ}, see also \cite{L2}).

\begin{definition}
Let $E$ be a Hermitian bundle over an $m$ dimensional K\"{a}hler manifold $C$
with complexified K\"{a}hler form $\omega^{\mathbb{C}}$. A differential form
$B\in\Omega^{0,q}\left(  C,ad\left(  E\right)  \right)  $ is called a deformed
$\bar{\partial}$-harmonic form if it satisfies
\begin{align*}
\bar{\partial}B  &  =0,\\
\operatorname{Im}e^{i\theta}\left(  \omega^{\mathbb{C}}+F\right)  ^{m-q}%
\wedge\partial B  &  =0.
\end{align*}
Here $F$ is the curvature tensor of $E$.

We denote the space of their solutions as $QH^{q}\left(  C,End\left(
E\right)  \right)  $.
\end{definition}

If the connection on $E$, the B-field on $C$ and the phase angle $\theta$ are
all trivial, then a deformed $\bar{\partial}$-harmonic form is just an
ordinary $\bar{\partial}$-harmonic form. It is useful to know if there is
always a unique deformed $\bar{\partial}$-harmonic representative for each
cohomology class in $H^{q}\left(  C,End\left(  E\right)  \right)  $. One might
need to require $\omega^{\mathbb{C}}+F$ to be sufficiently positive to ensure
ellipticity of the equation.

\bigskip

\textbf{Deformation of B-cycles }

The situation becomes more interesting when we deform both the complex
submanifold $C\subset M$ and a holomorphic vector bundle $E$ over $C$. The
curvature tensor of $E$ determines a cohomology class in $H^{1,1}\left(
C,End\left(  E\right)  \right)  $ which is called the Atiyah class of $E$.
Using the isomorphism
\[
H^{1,1}\left(  C,End\left(  E\right)  \right)  =Ext_{O_{C}}^{1}\left(
T_{C},End\left(  E\right)  \right)  ,
\]
we obtain an extension bundle $\mathbb{E,}$%
\[
0\rightarrow End\left(  E\right)  \rightarrow\mathbb{E}\rightarrow
T_{C}\rightarrow0\text{.}%
\]
The cohomology group $H^{1}\left(  C,\mathbb{E}\right)  $ is the space of
infinitesimal deformations of the pair $\left(  C,E\right)  $ with $C$ a
complex manifold and $E$ a holomorphic vector bundle over $C$. It fits into
the long exact sequence
\[
H^{0}\left(  C,T_{C}\right)  \rightarrow H^{1}\left(  C,End\left(  E\right)
\right)  \rightarrow H^{1}\left(  C,\mathbb{E}\right)  \rightarrow
H^{1}\left(  C,T_{C}\right)  ,
\]
which relates various deformation problems.

If we forget $E$ and look at abstract deformations of $C$ and deformations of
$C$ inside $M$. Their relationships are explained by an exact sequence
\[
H^{0}\left(  C,T_{M}|_{C}\right)  \rightarrow H^{0}\left(  C,N_{C/M}\right)
\rightarrow H^{1}\left(  C,T_{C}\right)  ,
\]
which is induced from the short exact sequence
\[
0\rightarrow T_{C}\rightarrow T_{M}|_{C}\rightarrow N_{C/M}\rightarrow0.
\]
These cohomology groups (i) $H^{0}\left(  C,T_{M}|_{C}\right)  $, (ii)
$H^{0}\left(  C,N_{C/M}\right)  $ and (iii) $H^{1}\left(  C,T_{C}\right)  $
parametrize infinitesimal deformations of (i) $C$ inside $M$ with fix complex
structure on $C$; (ii) $C$ inside $M$ with various complex structures on $C$
and (iii) complex structures on $C$.

Now the infinitesimal deformations of a pair $\left(  C,E\right)  $ with $C$ a
submanifold of $M$ and $E$ a holomorphic bundle over $C$ are parametrized by
the Cartesian product of $H^{1}\left(  C,\mathbb{E}\right)  $ and
$H^{0}\left(  C,N_{C/M}\right)  $ over $H^{1}\left(  C,T_{C}\right)  $:
\[
H^{1}\left(  C,\mathbb{E}\right)  \times_{H^{1}\left(  C,T_{C}\right)  }%
H^{0}\left(  C,N_{C/M}\right)
\]
It is pointed out by Thomas that this space is the same as $Ext^{1}\left(
\iota_{\ast}E,\iota_{\ast}E\right)  $ where $\iota:C\rightarrow M$ is the
inclusion morphism.

The deRham representative of an element in this space is given by $B\in
\Omega^{0,1}\left(  C,End\left(  E\right)  \right)  $ and $v\in\Omega
^{0}\left(  N_{C/M}\right)  $ satisfying (i) $\bar{\partial}v=0$ and
$\bar{\partial}B+\left(  \bar{\partial}v\right)  \lrcorner F=0$ inside
$\Omega^{0,2}\left(  C,End\left(  E\right)  \right)  $ where $F$ is the
curvature tensor of $E$.

Using Hodge decomposition we can choose harmonic representative in each class.
Namely we can find unique pair $\left(  B,v\right)  $ in any fixed class in
$H^{1}\left(  \mathbb{E}\right)  \times_{H^{1}\left(  T_{C}\right)  }$
$H^{0}\left(  N_{C/M}\right)  $ which satisfies $\bar{\partial}^{\ast}B=0$ and
$\bar{\partial}v=0$. We can also rewrite these two equations as $\partial
B\wedge\omega^{n-1}=0$ and $\bar{\partial}v=0$. Here $\partial$ denote the
$\left(  1,0\right)  $ component of the covariant derivative for $End\left(
E\right)  $. Recall that there can be obstructions for such infinitesimal
deformations to come from a honest family.

Finally we can write down the equations whose solutions parametrize
infinitesimal variations of B-cycles $\left(  C,E\right)  ,$ where the
Hermitian metric on $E$ obeys the deformed Hermitian-Yang-Mills equations.
They are a deformation of the above harmonic equations for $B\in\Omega
^{0,1}\left(  C,End\left(  E\right)  \right)  $ and $v\in\Omega^{0}\left(
N_{C/M}\right)  $:
\[
\left\{
\begin{array}
[c]{r}%
\bar{\partial}B+\left(  \bar{\partial}v\right)  \lrcorner F=0,\\
\left[  \partial B\left(  \omega^{\mathbb{C}}+F\right)  ^{m-1}\right]
_{sym}=0,\\
\bar{\partial}v=0.
\end{array}
\right.
\]

Remark: It is more natural to consider deformation of B-cycles inside $M$
together with the deformation of the complex structure on $M$. These
infinitesimal deformations are parametrized by the Cartesian product,%

\[
H^{1}\left(  C,\mathbb{E}\right)  \times_{H^{1}\left(  C,T_{C}\right)  }%
H^{1}\left(  T_{M}\rightarrow N_{C/M}\right)  ,
\]
if we ignore the deformed Hermitian-Yang-Mills equation.

\bigskip

\textbf{Moduli space of B-cycles and its correlation function}

Let us denote the moduli space of B-cycles on $M$ by $_{B}\mathcal{M}\left(
M\right)  $. As a moduli space of holomorphic objects, it carries a natural
complex structure. It also have a natural pre-symplectic form $_{B}\omega$. At
a given point $\left(  C,E\right)  $, we denote the curvature of $E$ as $F$.
If $\left(  v_{1},B_{1}\right)  $ and $\left(  v_{2},B_{2}\right)  $ are two
tangent vectors of $_{B}\mathcal{M}\left(  M\right)  $ at $\left(  C,E\right)
$, then their pairing with this pre-symplectic form is given by
\begin{align*}
&  _{B}\omega\left(  \left(  v_{1},B_{1}\right)  ,\left(  v_{2},B_{2}\right)
\right) \\
&  =\operatorname{Im}e^{i\theta}\int_{C}Tr\left[  \left(  \omega^{\mathbb{C}%
}+F\right)  ^{m-1}B_{1}B_{2}\right]  _{symm}+\int_{C}\omega\left(  v_{1}%
,v_{2}\right)  \omega^{m}\text{.}%
\end{align*}
If we replace $\omega$ by a large multiple and set the phase angle to zero,
then the dominating term does give a symplectic two form:
\[
\int_{C}TrB_{1}B_{2}\omega^{m-1}+\int_{C}\omega\left(  v_{1},v_{2}\right)
\omega^{m}\text{.}%
\]
on the moduli space of usual Hermitian-Yang-Mills bundles over $C$.

Next we assume $M$ is a Calabi-Yau manifold, its holomorphic volume form
$\Omega$ induces an $n$-form $_{B}\Omega$ on the moduli space $_{B}%
\mathcal{M}\left(  M\right)  $. In string theory terminology this is the
correlation function of the theory.

Let us first describe $_{B}\Omega$ in explicit terms. If $\left(  C,E\right)
$ is a B-cycle on $M$ of dimension $m$. Let
\[
\left(  B_{j},v_{j}\right)  \in\Omega^{0,1}\left(  C,End\left(  E\right)
\right)  \times\Omega^{0}\left(  N_{C/M}\right)
\]
with $j=1,2,...,n$ represent tangent vectors of $_{B}\mathcal{M}\left(
M\right)  $ at $\left(  C,E\right)  $ as before. Then we have%

\begin{align*}
&  _{B}\Omega\left(  C,E\right)  \left(  \left(  B_{j},v_{j}\right)  \right)
_{j=1}^{n}\\
&  =\sum_{\left(  j_{1},\ldots,j_{n}\right)  \in S_{n}}\int_{C}Tr\left[
B_{j_{1}}\cdots B_{j_{m}}\right]  _{sym}\iota_{v_{j_{m+1}}}\cdots
\iota_{v_{j_{n}}}\Omega.
\end{align*}
We can also write down the correlation function in an intrinsic way: Consider
the diagram
\[%
\begin{array}
[c]{ccc}%
C\times Map\left(  C,M\right)  \times\mathcal{A}_{C}\left(  E\right)  &
\overset{ev}{\rightarrow} & M\\
\downarrow\pi &  & \\
Map\left(  C,M\right)  \times\mathcal{A}_{C}\left(  E\right)  &  &
\end{array}
\]
where $\mathcal{A}_{C}\left(  E\right)  $ is the space of Hermitian
connections on $E$. Over $C\times\mathcal{A}_{C}\left(  E\right)  $ there is a
universal connection and we denote its curvature tensor as $\mathbb{F}$. There
is also an evaluation map $ev:C\times Map\left(  C,M\right)  \rightarrow M$.
Then the correlation function $_{B}\Omega$ can be expressed as follows (see
for example \cite{L2}),
\[
_{B}\Omega=\int_{C}Tr\mathbb{F}^{m}\wedge ev^{\ast}\Omega.
\]
This can be expressed simply as the composition,
\[
\otimes^{n}Ext^{1}\left(  \iota_{\ast}E,\iota_{\ast}E\right)  \rightarrow
Ext^{n}\left(  \iota_{\ast}E,\iota_{\ast}E\right)  \overset{Tr}{\rightarrow
}H^{n}\left(  M,O_{M}\right)  \cong\mathbb{C}\text{,}%
\]
provided that we can replace the deformed Hermitian-Yang-Mills equation by the
usual Hermitian-Yang-Mills equation.

\bigskip

\bigskip

\textbf{Symplectic reduction}

On the configuration space $Map\left(  C,M\right)  \times\mathcal{A}%
_{C}\left(  E\right)  $ we can write down a similar pre-symplectic form in an
intrinsic way:
\begin{align*}
_{B}\omega &  =\operatorname{Im}e^{i\theta}\left[  \int_{C}Tr\left(
\omega^{\mathbb{C}}+\mathbb{F}\right)  ^{m+1}\right]  ^{(2)}+\left[  \int
_{C}ev^{\ast}\omega^{m+1}\right]  ^{(2)}\\
&  =\frac{1}{\left(  m+1\right)  !}\left[  \int_{C}\operatorname{Im}%
Tre^{i\theta+\omega^{\mathbb{C}}+F+ev^{\ast}\omega}\right]  ^{(2)}\text{.}%
\end{align*}
The group of gauge transformations of $E$ acts on $\mathcal{A}_{C}\left(
E\right)  $, and therefore on $Map\left(  C,M\right)  \times\mathcal{A}%
_{C}\left(  E\right)  $ preserving $_{B}\omega$. The complex submanifold of
$Map\left(  C,M\right)  \times\mathcal{A}_{C}\left(  E\right)  $, consisting
of complex submanifolds $C\subset M$ and holomorphic connections (i.e.
$F^{2,0}=0$) on $E$, is invariant under this gauge group action.

The moment map equation is the deformed Hermitian-Yang-Mills equation and the
symplectic quotient is simply the moduli space of B-cycles $_{B}%
\mathcal{M}\left(  M\right)  $.

\bigskip

Email: leung@math.umn.edu
\end{document}